\documentclass[11pt,leqno]{article}
\usepackage{amsmath, amscd, amsthm, amssymb, graphics, xypic, mathrsfs,textcomp,url}
\usepackage{fancyhdr}
\usepackage[pagebackref=true]{hyperref}
\hypersetup{backref}

\newcommand{\setof}[1]{\{ #1 \}}

\newcommand{\aone}{{\mathbb A}^1}
\newcommand{\pone}{{\mathbb P}^1}

\newcommand{\tensor}{\otimes}

\renewcommand{\O}{{\mathcal O}}

\newcommand{\F}{{\mathcal F}}

\newcommand{\dmeff}[1]{{DM}^{eff,-}_{Nis}(k,{#1})}
\newcommand{\Shv}{{\mathcal Shv}}

\newcommand{\cplx}{{\mathbb C}}

\newcommand{\Z}{{\mathbb Z}}

\newcommand{\Q}{{\mathbb Q}}

\newcommand{\Sym}{{\mathsf{Sym}}}

\newcommand{\ga}{{\mathbb G}_{\bf a}}
\newcommand{\gm}{{\mathbb G}_{\bf m}}

\newcommand{\simpnis}{{\Delta}^{\circ}Shv_{Nis}({\mathcal Sm}/k)}

\newcommand{\isomto}{\stackrel{\sim}{\longrightarrow}}

\newcommand{\Top}{{\mathcal Top}}

\newcommand{\Spec}{\operatorname{Spec}}

\newcommand{\Sm}{{\mathcal Sm}}
\newcommand{\Sch}{{\mathcal Sch}}

\newcommand{\kbar}{{\overline{k}}}

\newcommand{\real}{{\mathbb R}}

\theoremstyle{plain}
\newtheorem{thm}{Theorem}[section]
\newtheorem{lem}[thm]{Lemma}
\newtheorem{keylem}[thm]{Key Lemma}
\newtheorem{cor}[thm]{Corollary}
\newtheorem{prop}[thm]{Proposition}
\newtheorem{claim}[thm]{Claim}
\newtheorem{question}[thm]{Question}

\newtheorem{conj}[thm]{Conjecture}

\theoremstyle{definition}
\newtheorem{defn}[thm]{Definition}

\theoremstyle{remark}
\newtheorem{rem}[thm]{Remark}
\newtheorem{ex}[thm]{Example}

\numberwithin{equation}{section}

\newcommand{\shrinkmargins}[1]{
  \addtolength{\textheight}{#1\topmargin}
  \addtolength{\textheight}{#1\topmargin}
  \addtolength{\textwidth}{#1\oddsidemargin}
  \addtolength{\textwidth}{#1\evensidemargin}
  \addtolength{\topmargin}{-#1\topmargin}
  \addtolength{\oddsidemargin}{-#1\oddsidemargin}
  \addtolength{\evensidemargin}{-#1\evensidemargin}
  }

\shrinkmargins{.9}

\begin{document}
\pagestyle{fancy}
\renewcommand{\sectionmark}[1]{\markright{\thesection\ #1}}
\fancyhead{}
\fancyhead[LO,RE]{\bfseries\footnotesize\thepage}
\fancyhead[LE]{\bfseries\footnotesize\rightmark}
\fancyhead[RO]{\bfseries\footnotesize\rightmark}
\chead[]{}
\cfoot[]{}
\setlength{\headheight}{1cm}

\title{{\bf On unipotent quotients and \\ some $\aone$-contractible smooth schemes}}
\date{}
\author{Aravind Asok \\ \begin{footnotesize}Department of Mathematics\end{footnotesize} \\ \begin{footnotesize}University of Washington\end{footnotesize} \\ \begin{footnotesize}Seattle, WA 98195 \end{footnotesize} \\ \begin{footnotesize}\url{asok@math.washington.edu}\end{footnotesize}\\
\and Brent Doran \\ \begin{footnotesize}School of
Mathematics\end{footnotesize} \\ \begin{footnotesize}Institute for
Advanced Study\end{footnotesize} \\ \begin{footnotesize}Princeton,
NJ 08540\end{footnotesize} \\
\begin{footnotesize}\url{doranb@ias.edu}\end{footnotesize}}
\maketitle

\begin{abstract}
We study quotients of quasi-affine schemes by unipotent groups over
fields of characteristic $0$.  To do this, we introduce a notion of
stability which allows us to characterize exactly when a principal
bundle quotient exists and, together with a cohomological vanishing
criterion, to characterize whether or not the resulting quasi-affine
quotient scheme is affine. We completely analyze the case of
$\ga$-invariant hypersurfaces in a linear $\ga$-representation $W$;
here the above characterizations admit simple geometric and
algebraic interpretations. As an application, we produce arbitrary
dimensional families of non-isomorphic smooth quasi-affine but not
affine $n$-dimensional varieties ($n \geq 6$) that are contractible
in the sense of $\aone$-homotopy theory.  Indeed, existence follows
without any computation; yet explicit defining equations for the
varieties depend only on knowing some linear $\ga$- and $SL_2$-
invariants, which, for a sufficiently large class, we provide.
Similarly, we produce infinitely many non-isomorphic examples in
dimensions $4$ and $5$.
 Over $\cplx$, the analytic spaces underlying these varieties are
non-isomorphic, non-Stein, topologically contractible and often
diffeomorphic to ${\mathbb C}^n$.
\end{abstract}

\begin{footnotesize}
\tableofcontents
\end{footnotesize}

\section{Introduction}
\subsubsection*{History and Motivation}
In 1935, J.H.C. Whitehead constructed, as a counterexample to his
``proof" of the $3$-dimensional Poincar\'e conjecture, the first
example of an open (i.e., non-compact and without boundary) contractible manifold
not homeomorphic to a ball (see \cite{Whitehead}). Subsequently,
D.R. McMillan produced infinitely many pairwise non-homeomorphic
open contractible smooth $3$-manifolds (see \cite{McMillan}).
Slightly earlier, Mazur and Poenaru had provided examples of
contractible open $4$-manifolds (Mazur's examples can be constructed
as smooth manifolds, see \cite{Mazur, Poenaru}).  Generalizing these
constructions, Curtis and Kwun  (see \cite{CurtisKwun}) showed that
there exist infinitely many pairwise non-homeomorphic, contractible,
open $n$-manifolds for every $n \geq 5$, and Glaser (see
\cite{Glaser}) showed that the same result held in dimension $4$.

Roughly contemporaneously, geometric topologists began to explore
the possibility of ``exotic" $PL$ and smooth structures compatible
with the usual topology on $\real^n$. Stallings proved (see
\cite{Stallings}) that if $M$ is an open contractible $n$-manifold
of dimension $n \geq 5$, simply connected at infinity\footnote{An
open manifold $M$ is said to be simply connected at infinity if for
every compact subset $C \subset M$, there is a compact subset $D$
such that $C \subset D \subset M$ and $M \backslash D$ is connected
and simply connected.}, then $M$ is $PL$-isomorphic to $\real^n$;
if further $M$ is smooth, then $M$ is in fact diffeomorphic to
$\real^n$ with its usual smooth structure.  (This result also
follows from the $h$-cobordism theorem if $n \geq 6$, see
\cite{Milnor} \S 9 Proposition A).  In other words, for $n \geq 5$,
$\real^n$ admits, up to the appropriate notion of isomorphism,
unique $PL$ and smooth structures.

Surprisingly, a simple shift in perspective allows one to construct,
at least in principle, all contractible manifolds; this will be the
motivating theme of this paper. It follows from results of McMillan,
Zeeman (see \cite{McMillanZeeman}) and Stallings (see
\cite{Stallings}) that {\em all} of the examples just discussed can
be realized as quotients of free $\real^k$ actions on $\real^{n+k}$
(for appropriate $n$ and $k$). In fact, Stallings (see {\em loc.
cit.} \S 4, \S 5 and Proposition 2.2) shows that any open
contractible (PL or smooth) $3$-manifold can be constructed as a
quotient of $\real^5$ by a free (PL or smooth) $\real^2$-action and,
more generally, for any $n \geq 4$ any open contractible (PL or
smooth) $n$-manifold can be constructed as a quotient of
$\real^{n+1}$ by a free (PL or smooth)
$\real$-action.\footnote{Strictly speaking, Stallings shows that the
product of the manifold with $\mathbb{R}^1$ (or $\mathbb{R}^2$ when
$n=3$) is PL-isomorphic or diffeomorphic to $\mathbb{R}^n$; however,
for contractible real manifolds, any $\mathbb{R}^k$-bundle is
trivial, so our statement is equivalent to his, and generalizes well
to our algebraic setting.}

The na\"{i}ve algebro-geometric analog of this question is whether
all smooth contractible complex algebraic varieties can be
constructed as quotients of $\mathbb{C}^n$ by the free algebraic
action of a unipotent group. All of the constructions just mentioned
are inherently ``topological" so it perhaps came as a shock that a
smooth contractible variety, apart from affine space, even exists.
The first example was given by Ramanujam in his landmark paper
\cite{Ramanujam}.\footnote{Ramanujam proved much more: any smooth,
complex algebraic surface whose underlying analytic space is
contractible and simply connected at infinity is necessarily {\em
algebraically} isomorphic to ${\mathbb A}^2$. Ramanujam's example,
being non-simply connected at infinity, was necessarily not
homeomorphic to $\real^4$.}  It, together with the fact that Zariski
cancellation holds in dimension $2$ (see \cite{Fu}), provides a
counter-example to this analog of the question.

Ramanujam's example was only the tip of the iceberg.  Other authors
showed that there exist many examples of contractible smooth
algebraic varieties in every (complex) dimension $\geq 2$ (see the
beautiful survey \cite{Za} for an overview and many more
references).  In this paper, we begin a study of contractible
algebraic varieties from the standpoint of motivic homotopy theory.
Rather, since topological contractibility only makes sense for
varieties defined over fields that are embeddable in the complex
numbers, we have to reformulate the notion of contractibility
appropriately.

Following Morel and Voevodsky (see \cite{MV}) we view the category
of smooth schemes as analogous to the category of topological spaces
with the affine line playing the role of the unit interval in
ordinary topology.  Morel and Voevodsky replace the category of
(e.g. locally contractible) topological spaces by the category of
(simplicial) Nisnevich sheaves on $\Sm/k$ (the category of smooth
manifolds is replaced by the Nisnevich sheaves corresponding to
smooth schemes), the notion of homeomorphism is replaced by
isomorphism of smooth schemes and finally, the usual topological
homotopy category is replaced by the Morel-Voevodsky
$\aone$-homotopy or ``motivic homotopy" category.  These analogies
are, of course, not perfect (as we shall explain), but hopefully
serve to guide intuition.

Our goal is to study $\aone$-contractible smooth algebraic
varieties, i.e. those varieties that are $\aone$-weakly equivalent
to $\Spec k$.  Essentially by construction ${\mathbb A}^n$ is
$\aone$-contractible.  However, we will see that there are many
examples of smooth algebraic varieties, not isomorphic to affine
space, that are $\aone$-contractible.  Henceforth, we call such a
variety an {\em exotic $\aone$-contractible variety}. We suggest
that this notion gives the ``correct" algebro-geometric analog of
our thematic question, namely:

\begin{question} Does every smooth $\mathbb{A}^1$-contractible
variety arise as a quotient of affine space by the free action of a
unipotent group?
\end{question}

First, however, one needs to produce examples of interesting
$\mathbb{A}^1$-contractible varieties.  We prove, in this spirit,
the following results.

\begin{thm}[See Theorem \ref{thm:dimension4}]
\label{thm:introfamilies} For every $m \geq 4$, there exists a
denumerably infinite collection of pairwise non-isomorphic $m$-dimensional exotic $\aone$-contractible varieties, each admitting an embedding into a
smooth affine variety with pure codimension $2$ smooth boundary.
\end{thm}
\begin{thm}[See Theorem \ref{thm:moduli}]
For every $m \geq 6$ and every $n \geq 0$:
\begin{itemize}
\item there exists a connected $n$-dimensional scheme $S$ and a smooth
morphism $f: X \longrightarrow S$ of relative dimension $m$, whose
fibers over $k$-points are $\aone$-contractible and quasi-affine,
not affine, and pairwise non-isomorphic.

\item The morphism $f: X \longrightarrow S$ admits a partial compactification to a flat family $\bar{f}:
\overline{X} \longrightarrow S$ whose fibers over $k$-points are
smooth affine varieties. Furthermore, for any $k$-point $t \in S$,
the map $X_t \rightarrow \bar{X}_t$ is an open immersion with a
smooth complement of codimension $\geq 2$.
\end{itemize}
\end{thm}

In other words, there exist arbitrary dimensional {\em moduli} of
$\aone$-contractible smooth varieties in dimension $\geq 6$. We
stress that these examples are completely explicit and
non-pathological. The families arise in a simple geometric manner:
as $\ga$-quotients of families of $\ga$-invariant hypersurfaces in a
fixed linear $\ga$-representation $W$.  The resulting quotient
varieties are complements of smooth codimension $2$ subvarieties in
smooth hypersurfaces in $\Spec k[W]^{\ga}$.  The simplest case, for
example, is the complement in an affine quadric four-fold (defined
by the vanishing of $x_1 x_4 - x_2 x_3 - x_5(x_5+1)$ in
$\mathbb{A}^5$) of an explicit embedded copy of $\mathbb{A}^2$
(defined by $x_1 = x_2 = 0, x_5 = -1$); see Remark
\ref{rem:Winkelex} for details.

When $k = \cplx$, we prove that it is impossible to construct exotic $\aone$-contractible varieties of dimension $\leq 2$ by our method (see Claims \ref{claim:dimension1} and \ref{claim:dimension2}).  Indeed, there exists a
unique up to isomorphism smooth $\aone$-contractible variety of
dimension $1$, namely $\aone$. It follows from results of several authors, that all the examples of $\aone$-contractible smooth surfaces we produce are necessarily isomorphic to the affine plane.  Therefore, only dimension $3$ seems mysterious.  In analogy with the topological setting, one may need to use explicit $(\ga)^2$ actions to study dimension $3$.

The motivic homotopy category of schemes over $\Spec \cplx$, admits
a ``topological realization functor" to the usual homotopy category
of topological spaces.  This realization functor takes $\aone$-weak
equivalences of smooth schemes to ordinary weak equivalences and, in
particular, the topological realization of an $\aone$-contractible
smooth variety is a contractible smooth manifold. We are unable to
produce examples of contractible algebraic varieties that are
provably not $\aone$-contractible.  Topological intuition encourages
us to believe that such varieties exist; however, ``motivic"
intuition related to the Hodge conjecture imposes very strong
topological restrictions on any such examples. Summarizing the above
discussion, we make the following conjecture.\footnote{Note added in
proof: we can now produce examples of smooth affine surfaces over
$\cplx$ which are topologically contractible but {\em not}
$\aone$-contractible.}

\begin{conj}
For every $m \geq 3$, and every $n \geq 0$, there exists a connected
$n$-dimensional scheme $S$ and a smooth morphism $f: X
\longrightarrow S$ of relative dimension $m$, whose fibers are
$\aone$-contractible and for a fixed field $k$, the fibers of $f$
over $k$-points of $S$ are all non-isomorphic.
\end{conj}

Finally, we claim that the topological characterization ``at
infinity" of the $PL$ or smooth structure of $\real^n$ for $n \geq
5$ gives rise to a natural question: can one give a motivic
topological characterization of affine space? As a first step in
this direction, one can try to define a notion of motivic homology at infinity.
One such notion was introduced by Wildeshaus in his paper ``Basic
properties of the boundary motive" (see \cite{Wild}).  All exotic
$\aone$-contractibles have, via Poincar\'e duality, motivic homology
at infinity (in Wildeshaus' sense) that of a motivic sphere of appropriate dimension (see Lemma \ref{lem:homatinfinity}).  A natural question to ask is whether there exists a good notion of an ``$\aone$-singular chain complex at infinity" and of an ``$\aone$-fundamental group at infinity," analogous to the usual singular chain complex at infinity or the fundamental group at infinity, that one might use to characterize when an $\aone$-contractible smooth variety is exotic.  In this direction, Morel (see \cite{MICM}) dreams of a ``motivic $s$-cobordism theorem:" a characterization of affine space as a smooth scheme should be a related consequence.

\subsubsection*{Contents}
As the techniques used in this paper have been introduced fairly
recently, we have endeavored to make the paper as self-contained as
possible.  We begin, in \S \ref{s:contractibility}, by making a brief review of
$\aone$-contractibility.  In particular, we state the main criterion
we use to check that a morphism is an $\aone$-weak equivalence (see
Lemma \ref{lem:extendedhominv}). In addition, we give the simplest
examples of (singular) $\aone$-contractible algebraic varieties. In
\S \ref{s:unipotents}, we discuss the relevant elements from geometric invariant
theory for non-reductive group actions as developed in \cite{DK}. In
particular, after reviewing some basic facts about unipotent groups,
we discuss a condition characterizing the existence of principal
bundle quotients by unipotent group actions (see Theorems
\ref{thm:stablequotient} and Theorem \ref{thm:affinequotient}).  To
keep this section self-contained, we have given complete proofs of
all the main results; the focus here is on quotients of everywhere
stable quasi-affine schemes and the material is essentially
orthogonal to that contained in \cite{DK}.

In \S \ref{s:hypersurfaces}, we study the simplest class of unipotent group actions:
$\ga$-actions.  The Jacobsen-Morozov theorem (see Theorem
\ref{thm:jacobsenmorozov}) essentially allows us to reduce the study
of $\ga$-actions to $SL_2$-actions.  We completely resolve the
question of when the (principal bundle) quotient of a
$\ga$-invariant hypersurface in a linear $\ga$-representation is an
affine variety, a strictly quasi-affine variety, or not even a
scheme. In particular, we give two explicit characterizations (see
Theorems \ref{thm:geomchar1} and \ref{thm:algchar}) which are used
to produce all the examples discussed in \S \ref{s:examples}. One curious
consequence is a natural decomposition of {\em any} $\ga$-invariant
function into a sum of an $SL_2$-invariant function and a
$\ga$-invariant function of a very particular sort; the authors are
not aware of a classical version of this statement in invariant
theory (see Theorem \ref{thm:algchar}).  In \S \ref{s:examples}, we prove the two
theorems stated in the introduction by using a very simple class of
$\ga$-equivariant linear embeddings of affine space (see Theorems
\ref{thm:dimension4} and \ref{thm:moduli}).  We also make a detailed
study of strictly quasi-affine quotients in small dimensions. In \S
\ref{s:conjectures}, we discuss various consequences of and conjectures related to the
notion of $\aone$-contractibility.  We emphasize here that the very
existence of the motivic homotopy category allows us to make very
strong statements about the motivic topology of exotic
$\aone$-contractible varieties.  In particular, we discuss briefly
the idea of motivic topology at infinity.  We close in the Appendix
(\S \ref{s:appendix}) with a summary of the main tools of the technique of
faithfully flat descent, its application to Borel transfer, and the
proof of the quite general and formal Theorem
\ref{thm:affinequotient}. Consequences of descent and its
applications are utilized throughout the paper; rather than
interrupt the main discussion with technical sidelights, we have
compiled the relevant facts there.

\subsubsection*{Conventions}
Throughout this paper the word ``field" will stand for ``field of
characteristic zero."  The word ``scheme" will mean separated scheme,
locally of finite type over a field $k$, the word ``variety" will
mean ``reduced scheme of finite type," and all group schemes will be linear algebraic $k$-groups.  Given a group scheme $G$ and a scheme $X$, we will say $X$ is a $G$-scheme if $X$ admits an algebraic left
$G$-action; that being said, the geometric quotient of $X$ by $G$,
if it exists as a scheme, will be written $X/G$ (rather than $G
\backslash X$).  The word ``free" applied to a $G$-action on a
scheme $X$ will always mean scheme-theoretically free $G$-action,
i.e. the action morphism $G \times X \longrightarrow X \times X$ is
a closed immersion.

If $G$ is a reductive group, and $X$ is a $G$-scheme which admits a categorical quotient by the $G$-action, then we denote this categorical quotient by $X//G$ following the convention due to Mumford in \cite{GIT}.  Given a group scheme $G$ and a scheme $X$, a $G$-torsor (sometimes called a principal $G$-bundle) on $X$ is a triple $({\mathscr P},\pi,G)$ consisting of a finite type, faithfully flat morphism of schemes $\pi: {\mathscr P} \longrightarrow X$, such that the canonical morphism $G \times {\mathscr P} \longrightarrow {\mathscr P} \times {\mathscr P}$ is an isomorphism onto ${\mathscr P} \times_X {\mathscr P}$.  Observe that with our conventions, in particular separatedness of schemes, it follows from, e.g., \cite{GIT} Lemma 0.6 that if $({\mathscr P},\pi,G)$ is a $G$-torsor, then $G$ acts freely on ${\mathscr P}$.

Finally, given a closed immersion group homomorphism of linear
algebraic groups $H \hookrightarrow G$ and an $H$-scheme $X$, we
write $G *_H X$ for the twisted or contracted product;  this is the
(algebraic space, see Remark \ref{rem:algebraicspace}) quotient of $G \times X$ by the free $H$-action
defined by $h \cdot (g,x) = (gh^{-1},h\cdot x)$.\footnote{We use this
notation in contrast to the common topological convention of writing
$G \times_H X$ in order to avoid possible confusion with fiber
products.  As amalgamated products are never considered in the
paper, we hope this notation will induce no confusion in the
reader.} If $X$ is quasi-affine we prove in the appendix (see
Corollary \ref{cor:sheafproperty} and Remarks \ref{rem:quasi-affine}
and \ref{rem:qaffmorphisms}) that this contracted product exists as
a scheme.

\subsubsection*{Acknowledgements}
The authors owe an intellectual debt to the examples of J. Winkelmann in his paper \cite{Wi}, which are the first strictly quasi-affine quotients of a $\mathbb{C}^n$ to appear in the literature; it was a pleasure, and a wonderful self-check, to see them arise so naturally from another perspective. Our study of $\aone$-contractibility came as a sidenote to discussion and collaboration with Frances Kirwan; this work can be thought of as living in the intersection of two concurrent projects with her, we
thank her for various discussions around this subject matter. We
would like to thank James Parson, Paul Hacking, S\'andor Kov\'acs,
Steve Mitchell, and Fabien Morel for helpful discussions and answering various questions we had. We also would like to thank William Stein for
providing facilities for computation.  Finally, we are grateful to the referee for his close reading of the text and thoughtful suggestions.

This material is based upon work supported by the National Science
Foundation under agreement No. DMS-0111298.  Any opinions, findings
and conclusions or recommendations expressed in this material are
those of the authors and do not necessarily reflect the views of the
National Science Foundation.

\section{Contractibility in topology and algebraic geometry}
\label{s:contractibility} Let us begin by reviewing some basic
notions of $\aone$-homotopy theory.  The general references for the
material in this section are \cite{MV} and \cite{VICM}.  Let
${\Sm}/k$ denote the category of finite type smooth schemes over $k$
and let ${\Sch}/k$ denote the category of all finite type schemes
over $k$.  We will only use the most basic definitions from
\cite{MV}: the main goals of this section are to explain Lemma
\ref{lem:extendedhominv}, to prove Lemma \ref{lem:topcontract} that
the topological realization of any $\aone$-weak equivalence of
smooth schemes over a field embeddable in $\cplx$ is actually a
topological weak equivalence, and to give a sense for how smooth
$\aone$-contractible schemes are different from the singular case.
The reader willing to take these results on faith can proceed
directly to the geometric constructions of the next sections.

\subsubsection*{Spaces}
The category $\Sm/k$ is not suitable for the purposes of homotopy
theory: for example, quotients by subspaces do not always exist in
this category.  Equip $\Sm/k$ with the Nisnevich topology and
consider the category $\Shv_{Nis}(\Sm/k)$; we refer to objects of
this category as {\em spaces}.  As the Nisnevich topology is
sub-canonical, every representable presheaf is a sheaf.  Therefore,
the Yoneda embedding $\Sm/k \longrightarrow \Shv_{Nis}(\Sm/k)$,
which sends a smooth scheme $X$ to the representable functor $U
\mapsto X(U)$, is fully faithful.  For a possibly non-smooth scheme
$X$, the functor $U \mapsto X(U)$ is also a Nisnevich sheaf and
gives a functor $\Sch/k \longrightarrow \Shv_{Nis}(\Sm/k)$; this
extended functor is not fully faithful as the following example
shows.

\begin{ex}
\label{ex:cusps}
Let $p$ and $q$ be coprime integers.  Let $X_{p,q}$ denote the Nisnevich sheaf attached to the cuspidal curve $x^p - y^q = 0 \subset {\mathbb A}^2$.  Normalization determines a morphism $\aone \longrightarrow X_{p,q}$;
we now show that the induced morphism of sheaves is an isomorphism.
In fact, the presheaves on $\Sm/k$ associated with these schemes are
isomorphic.

The map of sheaves $\aone \longrightarrow X_{p,q}$ is surjective  To see this, observe that the normalization of $X_{p,q}$ is $\aone$, and any morphism from a connected test scheme $T$ to $X_{p,q}$ factors through the normalization.  Indeed, any dominant morphism factors through the normalization by the universal property, and any non-dominant morphism has image a point and can thus be lifted.  To show the map of sheaves $\aone \longrightarrow X_{p,q}$ is injective, we have to show that the diagonal closed immersion $\Delta: \aone \longrightarrow \aone \times_{X_{p,q}} \aone$ is surjective and hence identifies $\aone$ with the underlying reduced scheme of the product.  This follows using the fact that the morphism $\aone \longrightarrow X_{p,q}$ is radiciel.  We thank James Parson for explaining this example to us.
\end{ex}

\subsubsection*{$\aone$-weak equivalences}
As in topology, it is often convenient to consider categories of
simplicial spaces; denote by $\simpnis$ the category of simplicial
Nisnevich sheaves.  The motivic homotopy category can be constructed
from either the category $\Shv_{Nis}(\Sm/k)$ or $\simpnis$ by
localizing at the class of $\aone$-weak equivalences.  For a precise
definition of $\aone$-weak equivalence see \cite{MV} \S 2.2
Definitions 2.1 and 2.2 and \S 2.3 Proposition 3.14.  Roughly
speaking, for every sheaf $X$, we invert the projection morphism $X
\times \aone \longrightarrow X$.  This localization forces many
additional morphisms to be weak equivalences.  Let us now give some
examples of $\aone$-weak equivalences.

We let $i_0: \Spec k \longrightarrow \aone$ and $i_1: \Spec k
\longrightarrow \aone$ denote the inclusion of $k$-rational points
$0$ and $1$.  Given a pair of spaces $X$ and $Y$ and two morphisms
$f,g: X \longrightarrow Y$, an {\em elementary $\aone$-homotopy}
from $f$ to $g$ is a morphism $H: X \times \aone \longrightarrow Y$
such that $H \circ i_0 = f$ and $H \circ i_1 = g$.   The morphisms
$f$ and $g$ are said to be {\em $\aone$-homotopic} if they can be
connected by a finite sequence of elementary $\aone$-homotopies.
Finally a morphism $f: X \longrightarrow Y$ is said be a {\em strict
$\aone$-homotopy equivalence} if there exists a morphism $g: Y
\longrightarrow X$ such that $f \circ g$ and $g \circ f$ are
$\aone$-homotopic to $Id_X$ and $Id_Y$.

\begin{lem}[\cite{MV} \S 2.3 Lemma 3.6]
\label{lem:stricthomotopy}
Any strict $\aone$-homotopy equivalence is an $\aone$-weak equivalence.
\end{lem}

\begin{defn}
\label{defn:aonecontractible} A space $X$ is said to be {\em
$\aone$-contractible} if the structure morphism $X \longrightarrow
\Spec k$ is an $\aone$-weak equivalence.  A smooth
$\aone$-contractible scheme not isomorphic to affine space will be
called an {\em exotic $\aone$-contractible} scheme.
\end{defn}

\begin{lem}[\cite{MV} \S 3.2 Example 2.3]
\label{lem:extendedhominv}
Let $\pi: Y \longrightarrow X$ be a Zariski locally trivial, smooth
morphism of smooth schemes with $\aone$-contractible fibers.  Then
$\pi$ is an $\aone$-weak equivalence.
\end{lem}

\subsubsection*{Comparison with topological contractibility}
We now explain the relation between $\aone$-contractibility and the
usual topological weak equivalences when $k$ is a field that admits
an embedding into $\cplx$.  We follow the discussion of \cite{DI};
let $\Top$ denote the category of {\em all} topological spaces with
continuous maps as morphisms.  The usual notion of open set endows
$\Top$ with the structure of a Grothendieck site.  We let ${\mathscr
H}$ denote the usual homotopy category of topological spaces.

Consider the site ${\Sm/\cplx}_{Nis}$.  Dugger shows (see \cite{DUH}
Proposition 8.1) how to construct a model category
$U(\Sm/\cplx_{Nis})_{\aone}$ that is Quillen equivalent to the
Morel-Voevodsky category (and hence the resulting homotopy categories
are isomorphic).  Given a smooth $\cplx$-scheme $X$, the assignment
$X \longrightarrow X(\cplx)$ sending $X$ to its set of complex
points equipped with the usual topology extends to an adjoint pair
of Quillen functors from $U(\Sm/\cplx_{Nis})_{\aone}$ to $\Top$.  In
particular, any such functor preserves weak equivalences of
cofibrant objects.  We let $t^{\cplx}$ denote the induced functor of
homotopy categories.  The next result then follows from the
fact that representable sheaves are cofibrant objects in
$U(\Sm/\cplx_{Nis})_{\aone}$.

\begin{lem}
\label{lem:topcontract}
If $f: X \longrightarrow Y$ is a morphism of smooth schemes which is
an $\aone$-weak equivalence, then the induced map $t^{\cplx}(f):
X(\cplx) \longrightarrow Y(\cplx)$ is a topological weak equivalence.
In particular, if $X$ is any $\aone$-contractible smooth scheme, then
$X(\cplx)$ is a contractible topological space.
\end{lem}

\subsubsection*{Contracting $\gm$-actions and singular $\aone$-contractible varieties}
A natural way to produce explicit strict $\aone$-homotopy
equivalences is to consider $\gm$-actions.  If $X$ is a variety
equipped with an algebraic $\gm$-action such that there is a unique
closed $\gm$-orbit that is a fixed-point, then one expects $X$ to be
contractible since, from the standpoint of Morse theory, all points
``flow toward the fixed points."  More generally, suppose $T$ is a
$k$-torus acting on a scheme $X$ with a unique closed $T$-orbit that
is furthermore a fixed point (which is necessarily a $k$-rational
point).  In this case, we will say that $X$ admits a {\em
contracting $T$-action}.

\begin{rem}
\label{lem:contraction} Suppose $T$ is a $k$-split torus.  If $X$ is
an affine $T$-scheme equipped with a contracting $T$-action, then
$X$ is an $\aone$-contractible scheme.

To see this, observe that the inclusion of the $T$-fixed point determines a $T$-equivariant morphism $\iota: \Spec k \longrightarrow X$ and the structure morphism $X \longrightarrow \Spec k$ (equivalently the categorical quotient morphism $X \longrightarrow \Spec k[X]^T$) is a
$T$-equivariant morphism as well.

Choose a ``generic" one-parameter subgroup $\mu: \gm \longrightarrow
T$ that has the same fixed-point locus as $T$ (we can do this
because $T$ is split).  Such a choice allows us to reduce to the case where $T = \gm$.  Consider the induced action morphism $\mu: \gm \times X \longrightarrow X$.  We claim that the action morphism $\mu: \gm \times
X\longrightarrow X$ extends to a morphism $\bar{\mu}: \aone \times X
\longrightarrow X$.  This follows from results of Hesselink on the existence of the concentrator scheme (see \cite{Hesselink} Defn. 4.2 and 5.5).

Finally, let us show that $\bar{\mu}$ defines an elementary
$\aone$-homotopy between the identity map $Id: X \longrightarrow X$
and the composite map $X \longrightarrow \Spec k
\stackrel{\iota}{\longrightarrow} X$.  The structure morphism $X
\longrightarrow \Spec k$ is an $\aone$-weak equivalence by Lemma
\ref{lem:stricthomotopy}, and therefore, by Definition
\ref{defn:aonecontractible}, $X$ is $\aone$-contractible.
\end{rem}

\begin{ex}
Consider the linear action of $\gm$ on ${\mathbb A}^n$ with weights
$a_1,\ldots,a_n$.  Assume further that each of the $a_i$ is a
strictly positive integer.  Then any closed $\gm$-stable subvariety
necessarily has the origin as the unique $\gm$-fixed point and hence
satisfies the hypotheses of Remark \ref{lem:contraction}.  It follows
from Example \ref{ex:cusps} that the case $n = 2$ only produces
examples isomorphic as spaces to affine space.
\end{ex}

\begin{ex}
By Luna's slice theorem (see \cite{Luslice} III.1 Corollaire 2),
every smooth affine $\gm$-variety with a contracting $\gm$-action is
necessarily $\gm$-equivariantly isomorphic to a vector space endowed
with a linear $\gm$-action.  Therefore, $\aone$-contractible {\em
smooth} varieties constructed by means of Remark \ref{lem:contraction} are scheme-theoretically isomorphic to affine space.
\end{ex}

\begin{ex}
Note that the assumption that $X$ is affine is necessary in the
statement of Remark \ref{lem:contraction} as the following example
shows.  Take $\pone$ with the usual $\gm$-action and consider the
quotient that identifies the two $\gm$-fixed points.  The resulting
variety is non-normal (though semi-normal), and admits a $\gm$-action that has a unique closed orbit that is a $\gm$-fixed point.  One can check that this scheme is not $\aone$-contractible.
\end{ex}

\section{Stability for unipotent groups}
\label{s:unipotents}
In this section, we investigate a general technique for studying
quasi-affine schemes that have the structure of a $U$-torsor over
some base (where $U$ is a unipotent $k$-group).  We use a notion of
stability for $U$-actions in the spirit of the geometric invariant
theory for reductive groups (see \ref{defn:stable}).  We begin by
recalling some general facts about unipotent groups.  Then we prove
the basic results about stability and in particular about
``everywhere stable actions" on quasi-affine schemes (see Definition
\ref{defn:everywherestable} and Theorem \ref{thm:stablequotient}).
We present, following Greuel-Pfister and Kambayashi-Miyanishi-Takeuchi, a characterization of when the quotient of an affine $U$-scheme is affine (see Theorem \ref{thm:affinequotient}) but defer the proof to the Appendix, as it uses different language than the body of the paper (see \S
\ref{ss:affinequotient}).  Taken together, these results give a way
to characterize when quotients of an affine variety are affine,
quasi-affine, or not even a scheme; Corollary \ref{cor:main}
summarizes the results and suggests how they apply to constructing
$\aone$-contractible schemes.

\subsubsection*{Unipotent groups and a key lemma}
Let $U$ be a connected unipotent group over $k$. Recall the
following structure theorem; this is the essential ingredient in
understanding the structure of $U$-torsors.

\begin{thm}[Lazard (see \cite{KMT} Theorem 8.0)]
\label{thm:flag}
Let $U$ be a connected unipotent $k$-group.  Then $U$ admits an increasing filtration $F_i(U)$ by closed subgroups with successive
quotients $F_i(U)/F_{i-1}(U) \cong \ga$.
\end{thm}

\begin{cor}[Grothendieck]
\label{cor:unipgrptriv}
Let $X$ be a $k$-scheme and let $U$ be a
connected unipotent $k$-group.  Then all $U$-torsors on $X$ are
Zariski locally trivial.
\end{cor}

We defer the proof of this result to the appendix (see \S\ref{ss:groththeorem}).

\begin{keylem}
\label{lem:contractibility} Let $U$ be a connected unipotent group.
Suppose $U$ acts freely on an $\aone$-contractible finite type smooth scheme $X$ such that a geometric quotient $\pi: X \longrightarrow X/U$ exists as a scheme.  Then $X/U$ is smooth and $\aone$-contractible.
\end{keylem}

\begin{proof}
It follows from \cite{GIT} Proposition 0.9, that in this situation
the triple $(X,\pi,U)$ is a $U$-torsor.  We know that $\pi$ is a
Zariski locally trivial morphism by Corollary \ref{cor:unipgrptriv}.
Under our assumptions, unipotent groups are isomorphic to affine spaces, the fibers of $\pi$ are thus isomorphic to affine spaces and are hence
$\aone$-contractible.  If $X$ and $X/U$ are smooth, then $\pi$ is an
$\aone$-weak equivalence by Lemma \ref{lem:extendedhominv}.  Since
$X$ is $\aone$-contractible by assumption, the result would follow.
One need only observe that given a $U$-torsor $X \longrightarrow
X/U$, with $X/U$ a scheme, $X$ is smooth if and only if $X/U$ is
smooth.
\end{proof}

\subsubsection*{Stability and $U$-torsors}
Next, we discuss how to construct principal bundle quotients of
actions of unipotent groups using elements of the geometric
invariant theory for non-reductive groups studied by the second
author and F. Kirwan in \cite{DK}. In particular, we will show that
quotients satisfying the hypotheses of Lemma
\ref{lem:contractibility} abound. Suppose $U$ is a connected
unipotent group.  Any such group can be realized as a closed subgroup of a reductive group $G$; in fact every unipotent group is isomorphic to a closed subgroup of $GL_n$ for $n$ sufficiently large.  Let $i: U \hookrightarrow G$ denote the corresponding closed immersion group homomorphism.  We can now define a notion
of stability for actions of unipotent groups using the corresponding notion for
reductive groups.

Let $X$ be a quasi-affine $U$-scheme.  Consider the contracted product scheme $G *_U X$ (which exists by Corollary \ref{cor:sheafproperty} and Remark \ref{rem:quasi-affine}) with structure sheaf $\O_{G *_U X}$.  In this scenario, there is a canonical closed immersion $\iota: X \hookrightarrow G *_U X$ which sends a point $x \in X$ to the point $[e,x]$.  If the character group
$X^*(G)$ is trivial (e.g., $G$ is semi-simple), then $\O_{G*_U X}$
actually admits a canonical $G$-linearization.  For us it is
convenient to always work with the $G$-linearization associated with
the trivial character. We can then make the following definition.

\begin{defn}
\label{defn:stable}
A geometric point $x \in X$ is {\em stable}, denoted $x \in X^{s}$,
if for every reductive group $G$, and every closed immersion group homomorphism $i: U \hookrightarrow G$, the
geometric point $[e,x] = \iota(x)$ is a (properly) stable point in
the sense of Mumford (see \cite{GIT}) of $G *_U X$ with respect to
the $G$-linearized sheaf $\O_{G *_U X}$, where the linearization is
given by the trivial character.
\end{defn}

\begin{rem}
In the above definition, we can replace the unipotent group $U$ by any linear algebraic group.  In addition, the set of stable points is the set of geometric points of an open subscheme of $X$; the notation $X^s$ refers to this open subscheme.  As noted in \cite{GIT} Proposition 1.14, if $K/k$ is a field extension, then we have an identification $X^s_{K} \cong X^s_k \times_{k} K$ so that stability is preserved under field extensions.  We will therefore often make arguments by passing to an algebraic closure and then applying Galois descent; note that this is a special case of the descent theory described in Theorem \ref{thm:ffdescent}.
\end{rem}

Our definition of stability is intrinsic to the $U$-action on $X$, but it is not
{\em a priori} obvious that the set of stable points on $X$ is ever
non-empty.  We now show that to check whether a point is stable, it
suffices to check stability for a {\em single} reductive group.

\begin{lem}[Doran-Kirwan (compare \cite{DK} Lemma 5.1.6)]
\label{lem:semistableequalsstable} Given a reductive group $G$, a closed embedding group homomorphism $i: U \hookrightarrow G$ and
a geometric point $x \in X$, the point $\iota(x) = [e,x]$ in $G *_U
X$ is stable with respect to the linearization corresponding to the
trivial character, if and only if it is semi-stable in the sense of
Mumford with respect to the same linearization; i.e. $(G *_U X)^{s}
= (G *_U X)^{ss}$.
\end{lem}

\begin{proof}
Let $O$ be a $U$-orbit of geometric points in $X$ such that $G *_U
O$ is strictly semi-stable, i.e., lies in the locally closed
subscheme $X^{ss} \setminus X^s$. Consequently $G *_U O$ is a
$G$-orbit in the complement of a $G$-stable affine hypersurface
defined by a $G$-invariant polynomial $F$ (i.e., $G *_U O \subset (G
*_U X)_F$), and furthermore is either closed of non-maximal
dimension or not closed in $(G *_U X)_F$.  In either case, there is
a unique closed orbit $(G *_U O')$ in the closure of $G *_U O$ in
$(G *_U X)_F$ that is of non maximal dimension.  Now, a closed
sub-scheme of an affine scheme is affine, and hence the $G$-orbit
$(G *_U O')$ is necessarily affine; this means any point $y \in G
*_U O'$ necessarily has  strictly positive dimensional reductive
stabilizer group.  However, the stabilizer group of $y$ must be
conjugate to a subgroup of $U$ and is therefore unipotent as well;
it follows that the stabilizer group is in fact trivial, which is a
contradiction. Thus such an $O$ does not exist.
\end{proof}

\begin{prop}[Doran-Kirwan (compare \cite{DK} Proposition 5.1.8)]
\label{prop:stable-for-red}
A geometric point $x \in X$ is {\em stable} if and only if, for any fixed reductive group $G$, together with a fixed closed embedding group homomorphism $i: U \hookrightarrow G$, the point $\iota(x) = [e,x]$ in $G *_U X$ is
(properly) stable in the sense of Mumford with respect to the
linearization on $\O_{G *_U X}$ corresponding to the trivial
character of $G$.
\end{prop}

\begin{proof}
By Lemma \ref{lem:semistableequalsstable} it suffices to show that
if $G_1$ and $G_2$ are two reductive groups both containing $U$ that
the intersections $\iota_1(X) \cap (G_1 *_U X)^{ss}$ and $\iota_2(X)
\cap (G_2 *_U X)^{ss}$ co-incide.  By definition, a point $y$ in
$G_i *_U X$ is semi-stable if and only if it is contained in a
$G_i$-invariant open affine $(G_i *_U X)_F$, for $F$ a $G_i$-invariant
function in $k[G_i *_U X]$.

Fix a reductive group $G$ containing both $G_1$ and $G_2$ (e.g. $G_1
\times G_2$).  Using the Borel transfer (see \S \ref{ss:transfer}), we have isomorphisms:
$$
k[G *_U X]^G \isomto k[G_i *_U X]^{G_i} \isomto k[X]^U.
$$
If $f$ is a $U$-invariant, we let $F_i$ be the corresponding
$G_i$-invariants in $k[G_i *_U X]$ and $F$ the corresponding
$G$-invariant in $k[G *_U X]$ obtained by the transfer isomorphisms.
It follows that $\iota_1(x) \in (G_1 *_U X)_{F_1}$ if and only if
$\iota_2(x) \in (G_2 *_U X)_{F_2}$.  It suffices then to check that
the hypersurface complement $(G_1 *_U X)_{F_1}$ is affine if and
only if $(G_2 *_U X)_{F_2}$. To do this, one need only check that
$(G_i *_{U} X)_{F_i}$ is affine if and only if $(G *_U X)_{F}$ is
affine.

Now, since $G *_U X = G *_{G_i} (G_i *_U X)$ we know that $G*_U X$
has the structure of an \'etale locally trivial fiber space over
$G/G_i$.  Therefore, using the transfer isomorphism we have an
identification $(G *_U X)_F = G *_{G_i} (G_i *_U X)_{F_i}$.  As
$G_i$ is reductive, $G/G_i$ is affine by Matsushima's theorem.  The
projection onto the first factor makes $G *_{G_i} (G_i *_U X)_{F_i}$
an \'etale locally trivial fiber bundle over the affine variety
$G/G_i$ with fibers isomorphic to $(G_i *_U X)_{F_i}$.  By Corollary \ref{cor:morphismproperty}, the result follows.
\end{proof}

\begin{defn}
\label{defn:everywherestable}
Suppose a linear algebraic group $U$ acts on a quasi-affine scheme $X$.
We say that the action is {\em everywhere stable} if every geometric
point of $X$ is stable, i.e., if $X = X^s$.
\end{defn}

\begin{rem}
\label{rem:algebraicspace}
Artin has shown (see \cite{Artin} 6.3 Corollary) that if an affine algebraic group $G$ acts freely on a scheme $X$ that a quotient of $X$ by $G$ exists as an algebraic space.  In particular, if a unipotent group acts everywhere stably on a quasi-affine scheme $X$, the proof of Theorem \ref{thm:stablequotient} will show that the action is in fact free.  We also give an example of a free action of a unipotent group on an affine scheme for which the geometric quotient is an algebraic space but not a scheme (in particular, the action is not everywhere stable, see \ref{ex:freealgspace}).
\end{rem}

\begin{thm}[Doran-Kirwan (compare \cite{DK} Theorem 5.3.1)]
\label{thm:stablequotient} Let $U$ be a connected unipotent group and suppose $X$ is a finite type quasi-affine $U$-scheme. The
action of $U$ on $X$ is everywhere stable if and only if the
quotient morphism $\pi: X \longrightarrow X/U$ is a $U$-torsor with $X/U$ a
quasi-affine scheme.  In this case, the scheme $X/U$ can be
realized as an open subscheme of the scheme $\Spec k[X]^U$.  If furthermore $X$ is a variety, then $X/U$ is a variety as well.
\end{thm}

\begin{proof}
First, note that any finite subgroup of a connected unipotent group over a
characteristic zero field is trivial. Since the $U$-action on $X$ is
everywhere stable, it is proper. Since  $U$ has no finite subgroups,
the $U$-action on $X$ is set-theoretically free as well.  By Lemma \ref{lem:freeproper} the $U$-action on $X$ is free.

Note that since $U$ acts everywhere stably on $X$, $G$ acts everywhere stably on $G *_U X$; by the discussion of the previous paragraph, $G$ must actually act freely on $G *_U X$.  Since $G$ acts everywhere stably on $G *_U X$, we know by \cite{GIT} Chapter 1 Theorem 1.10 that a geometric quotient $G *_U X/G$ exists as a scheme and because the $G$-action is free, it follows that the quotient morphism $q: G *_U X \longrightarrow G *_U X/G$ is a $G$-torsor.

We begin with the forward implication.  We claim that $G *_U X/G$ is in fact a geometric quotient of $X$ by $U$, so that we have $G *_U X/G \isomto X/U$ as schemes.

Since $G *_U X/G$ is a geometric quotient by $G$,
it follows that $k$-points in $G *_U X/G$ correspond bijectively to
$G$-orbits in $G *_U X$.  By faithfully flat descent (see \ref{cor:morphismproperty}), closed
$G$-stable subschemes of $G *_U X$ are in bijection with closed
$U$-stable subschemes of $X$.  Therefore, $k$-points of $G *_U X/G$
are in bijection with $U$-orbits in $X$.  Next, $q: G *_U X \longrightarrow G*_U X/G$ is a submersive morphism hence a subscheme $V
\subset G *_U X/G$ is open if and only if $q^{-1}(V)$ is open in $G
*_U X$. Now, $q^{-1}(V)$ is a $G$-stable open and hence corresponds
to a $U$-stable open of $X$ (again by Corollary \ref{cor:morphismproperty}), and agrees with $\pi^{-1}(V)$ by definition.  Finally if $\iota: X \longrightarrow G*_U X$, we note by Corollary \ref{cor:sheafproperty} combined with \cite{SGA1} Expose VIII Cor. 1.9, that $q_*(\O_{G *_U X})^G \cong \pi_*(\O_{X})^U$.  Therefore, $G *_U X/G$ is in fact a geometric quotient of $X$ by $U$, and we will write $X/U$ and $G *_U X/G$ interchangeably henceforth.

Now we observe that $X \longrightarrow X/U$ is in fact a $U$-torsor. Indeed, $G *_U X \longrightarrow X/U$ is a $G$-torsor
over $X/U$. By pull-back the scheme $G \times X$ is a $G \times
U$-torsor over $X/U$. Therefore, by faithfully flat descent, $X$ is
a $U$-torsor over $X/U$.  (In fact, the same argument shows that $X
\longrightarrow X/U$ is a $U$-torsor if and only if $G *_U X
\longrightarrow G *_U X/G$ is a $G$-torsor.)  Since $X \longrightarrow X/U$ is a $U$-torsor it is necessarily an affine morphism.

Let us now check that $X/U$ is in fact a quasi-affine scheme.  To do this we check that $G *_U X$ embeds $G$-equivariantly in an affine scheme $\overline{G *_U X}$ such that $G *_U X/G \hookrightarrow \overline{G *_U X}/\!/G$.  We can pick a finite set $S$ of $G$-invariant functions in $k[G *_U X]$ such that given geometric points $x, y \in G *_U X$ if there exists some $G$-invariant function $f$ such $f(x) \neq f(y)$ then in fact we can find $f' \in S$ such that $f'(x) \neq f'(y)$; this follows from the fact that $G *_U X$ is of finite type (see Corollary \ref{cor:morphismproperty}).  We then consider the morphism $G *_U X \longrightarrow {\mathbb A}^k$ defined by the functions $f_i \in S$.  By definition, this gives an embedding $G *_U X/G$ into $\Spec k[f_1,\ldots,f_n]$.  We denote this closure by $\overline{G *_U X}/\!/G$, and this identifies $X/U$ as a quasi-affine scheme.

Now, for the reverse direction, suppose $X \longrightarrow X/U$ is a
principal bundle quotient.   For any $G$-invariant function $F$, the
induced morphism $(G *_U X)_F \longrightarrow (G *_U X)_F/G$ is
necessarily a principal $G$-bundle, and if $f$ denotes the
$U$-invariant in $k[X]$ corresponding to $F$, $X_f \longrightarrow
X_f/U$ is a $U$-torsor. If $X_f/U$ is affine, then $(G *_U X)_F/G$
is necessarily affine. Therefore, if $x$ is a geometric point of
$X_f$, $\iota(x) = [e,x]$ is a geometric point of  $(G *_U X)_F$.
Therefore, $\iota(x) = [e,x]$ is semi-stable in the sense of
Mumford. The result follows then by Lemma
\ref{lem:semistableequalsstable}.

Zariski's Main Theorem shows that we have an open immersion $i_G: G
*_U X/G \hookrightarrow \Spec k[G *_U X/G]$. By properties of a
geometric quotient, $\Spec k[G *_U X/G] \isomto \Spec k[G *_U X]^G$;
we thus obtain a morphism $X/U \hookrightarrow \Spec k[X]^U$, which
is {\em a fortiori} an open immersion because $i_G$ has the same
property. It follows that the scheme theoretic image of $X$ under
$\pi$ is a variety which is an open subscheme of $\Spec k[X]^U$.
\end{proof}

\begin{lem}
\label{lem:freeproper}
If $G$ is a linear algebraic group acting properly and set-theoretically freely on a scheme $X$, then the action is free.
\end{lem}

\begin{proof}
We must show that the action map $G \times X \longrightarrow X \times X$ is a closed embedding.  By definition of properness of an action, the action map is proper and quasi-finite and thus finite.  Since the stabilizers are trivial, it is unramified so it is an embedding.  Finally, the map is injective on geometric points and hence an embedding.
\end{proof}

\begin{rem}
Observe that Theorem \ref{thm:stablequotient} implies that the quotient of an everywhere stable action of a unipotent group $U$ on a strictly quasi-affine scheme $X$ (i.e. quasi-affine but not affine) is necessarily strictly quasi-affine.  Indeed, the morphism $X \longrightarrow X/U$ is affine and so if $X/U$ were affine, since the composite of two affine morphisms is affine, this would mean that $X$ was affine.
\end{rem}

\begin{rem}
Assume that $X$ is an $\aone$-contractible smooth affine scheme
(e.g. ${\mathbb A}^n$).  By Lemma \ref{lem:contractibility}, if $X$
admits an everywhere stable action of a unipotent group, the
quotient $X/U$ exists as a scheme and is smooth and
$\aone$-contractible.  If the quotient $X/U$ is {\em not} affine,
then it cannot be isomorphic to affine space, and so is an exotic
$\aone$-contractible scheme.
\end{rem}

Combining the above remarks, if $X$ is an affine scheme equipped with an everywhere stable $U$-action, we would like a criterion which characterizes when the quotient is an affine scheme.  There is an effective
cohomological criterion for this, which we adapt from a theorem of Greuel and Pfister (see \cite{GP} Theorem 3.10) or earlier by Kambayashi, Miyanishi and Takeuchi (see \cite{KMT} Theorem 7.1.1).  To state the result, we need some notation.

Let ${\mathfrak u}$ denote the Lie algebra corresponding to $U$; it
is necessarily nilpotent. The exponential map then defines an {\em
algebraic} isomorphism $\exp: \mathfrak{u} \isomto U$.  Now,
specifying a $U$-action on an affine $k$-scheme $X$ is equivalent to
giving a map $U \hookrightarrow Aut_k(X)$. Such an action determines
a Lie algebra homomorphism ${\mathfrak u} \longrightarrow
Der_k(k[X])$. The image of this map necessarily consists of locally
nilpotent derivations. Conversely, a unipotent group action can be
specified completely by giving a set of locally nilpotent
derivations generating the ${\mathfrak u}$-action.  The following
theorem is extremely useful, but we defer the proof to the appendix
(see \S \ref{ss:affinequotient}).

\begin{thm}
\label{thm:affinequotient} Suppose $X$ is an affine scheme equipped
with an action of a connected unipotent group $U$.  The following conditions are equivalent:
\begin{itemize}
\item[i)] The quotient map $\pi: X \longrightarrow X/U$ is a $U$-torsor and $\pi$ induces an isomorphism $X/U \cong \Spec(k[X]^U)$, i.e. $X/U$ is an affine scheme.
\item[ii)] The Lie algebra cohomology $H^1({\mathfrak u}, \Gamma(X,\O_X)) = 0$.
\item[iii)] The quotient map $X \longrightarrow X/U$ is a trivial $U$-torsor over the affine scheme $X/U$.
\end{itemize}
\end{thm}

\begin{ex}
\label{ex:imderivation}
Consider the situation where $X$ is an affine $\ga$-scheme.  Then the condition $H^1({\mathfrak g}_a,k[X]) = 0$ can be made more explicit.  Let $D$ be the locally nilpotent derivation defining the action of ${\mathfrak g}_a$ on $k[X]$.  Then, the ${\mathfrak g}_a$-cohomology of $k[X]$ can be computed as the cohomology of the Chevalley complex:
$$
Hom(k,k[X]) \stackrel{\delta}{\longrightarrow} Hom({\mathfrak g}_a,k[X])
$$
where the differential $\delta$ is defined by $\delta(\omega)(f) = \omega(Df)$.  Because $H^1({\mathfrak g}_a,k[X]) = Coker(\delta)$, we know that $H^1({\mathfrak g}_a,k[X]) = 0$ only if $\delta$ is surjective.  Therefore, viewing $D$ as a map $k[X] \longrightarrow k[X]$, $H^1({\mathfrak g}_a,k[X])$ is zero if and only if $1 \in Im(D)$.  If $1 \in Im(D)$, an element $s \in k[X]$ such that $D(s) = 1$ is sometimes called a {\em slice}.  Assuming $k[X]^{\ga}$ is finitely generated, one can give an algorithm to determine whether or not $1 \in Im(D)$.
\end{ex}

\subsubsection*{Spaces of actions}
Before we proceed, let us summarize the state of unipotent group
quotients as it now stands for us.  Fix an affine scheme $X$ and a
unipotent group $U$ whose associated Lie algebra is ${\mathfrak u}$.
Consider the set of locally nilpotent derivations
$Der_k^{ln}(k[X])$; this set is not in general a $k$-vector space.
Nevertheless, specifying a $U$-action on $X$ is equivalent to
specifying a Lie sub-algebra of $Der_k(k[X])$ isomorphic to
${\mathfrak u}$ that consists of locally nilpotent derivations.  We
will pay special attention to the case $U = \ga$, where the set of
actions of $\ga$ on $X$ is equivalent to the set $Der_k^{ln}(k[X])$.
We have the following schematic.
\begin{eqnarray*} H^1({\mathfrak u}, k[X]) = 0 \text{ actions} \subset \text{everywhere stable
actions}
\subset \text{free actions} \subset \text{actions}
\end{eqnarray*}
Thus, the inclusions above correspond to quotients which are:
\begin{eqnarray*} \text{affine varieties (geometric quotient)} \subset
\text{quasi-affine varieties (geometric quotient)} \\ \subset
\text{algebraic spaces
(geometric quotient)}  \subset \text{categorical quotient needn't exist}
\end{eqnarray*}

What is more, each inclusion above is strict, in that we can give
examples in each class not lying in the previous class.  Indeed, in principle there are ``algorithms" to detect which case holds.  The rest of the paper will be devoted to studying the first two inclusions, so let us give examples of the
other two; indeed, to show how widespread these phenomena are, we
give examples for $\ga$-actions on affine space $\mathbb{A}^5$.

\begin{ex}
\label{ex:freealgspace}
We now construct a free action of $\ga$ on an affine variety, in
fact $\mathbb{A}^5$, for which a quotient exists only as an
algebraic space. The rough idea, which we do not justify here, is to
find a $\ga$-invariant subvariety $X$ of a linear
$\ga$-representation $W$ such that (i) the action is free, (ii) it
has an open subset of stable points and (iii) it has a positive
dimensional subset of unstable points (i.e., points where all
homogeneous invariants vanish).  All of these can be explicitly
checked.  Here is one such construction, which recovers the example
due to Deveney and Finston \cite{DF}.

Let $V$ be the 2-dimensional representation of $\ga$ defined by the
usual embedding of $\ga$ into $SL_2$ as strictly lower triangular
matrices with ones along the diagonal.  Let $W = V \oplus V \oplus
\Sym^3 V \cong \mathbb{A}^8$ with coordinates $w_1, \ldots w_8$. Let
$X \cong \mathbb{A}^5$ be presented as a codimension $3$ closed
subvariety with defining equations $w_5 = 2w_1 w_3^2, w_6 = 2 w_1
w_3 w_4$, and $w_7 = 1 + w_1 w_4^2$.  One can check the action on
$X$ is proper and set-theoretically free; hence, it is free. It has
a $3$-dimensional subspace, determined by $w_1 = w_3 = 0$, in the
``unstable" locus.
\end{ex}

\begin{ex}
Let $V$ be the 2-dimensional representation of $\ga$ defined by the
usual embedding of $\ga$ into $SL_2$ as strictly lower triangular
matrices with ones along the diagonal.  Consider the $\ga$-action on
$\Sym^4 V \cong \mathbb{A}^5$. Then a categorical quotient of
$\Sym^4 V$ by the $\ga$-action does not exist.   In particular,
there are not enough invariants to separate any closed orbits
occurring in certain families. So the ``quotient" with respect to
invariants turns out to just be a constructible set, albeit sitting
in a slightly larger canonical quasi-affine variety.  If a
categorical quotient $Y$ of $X$ by a group $U$ exists, then the
quotient morphism $X \rightarrow Y$ must be surjective; thus here a
categorical quotient can not exist. See \cite{DK} Section $6$ for
details (in the projective case, although the arguments are the
same).  Note that if $G$ is a reductive group and $X$ is an affine
$G$-scheme, then while similar ``dimensional collapsing" phenomena
still occur, the morphism $X \rightarrow \Spec k[X]^G$ is always
surjective, so this problem involving constructible sets never occurs.
\end{ex}

In summary, the following corollary of the above theorems is computationally effective:

\begin{cor}
\label{cor:main}
Suppose a unipotent group $U$ acts freely on a finite type smooth affine scheme $X$.
\begin{itemize}
\item[i)] The algebraic space quotient $X/U$ is a smooth quasi-affine scheme if and only if the action of $U$ on $X$ is {\em everywhere stable}.
\item[ii)] The algebraic space $X/U$ is an affine scheme if and only if $H^1({\mathfrak u},k[X]) = 0$ in which case $k[X/U] = \Spec k[X]^U$, i.e. $X \longrightarrow \Spec k[X]^U$ is surjective; and $k[X]^U$ is necessarily finitely generated.
\item[iii)] Assume further that the action of $U$ on $X$ is everywhere stable.  If $H^1({\mathfrak u},k[X])$ is non-zero, then the smooth scheme $X/U$ has complement of codimension $\geq 2$ in $\Spec k[X]^U$.
\item[iv)] Given a unipotent group $U$ and two everywhere stable actions of $U$ on $X$, the resulting quotients $X/U$ are $\aone$-weakly equivalent.  If $X$ is an $\aone$-contractible variety, then $X/U$ is an $\aone$-contractible variety.  If $H^1({\mathfrak u},k[X])$ is non-trivial, then $X/U$ is not isomorphic to affine space.
\end{itemize}
\end{cor}

\begin{proof}
The only statement that doesn't follow immediately from Theorems
\ref{thm:stablequotient}, \ref{thm:affinequotient} and Lemmas
\ref{lem:contractibility} and \ref{lem:extendedhominv} is (iii).  Again by Theorems
\ref{thm:affinequotient} (i) and \ref{thm:stablequotient} (i) we can
assume that the morphism $j: X/U \longrightarrow \Spec k[X]^U$ is an
open embedding and not an isomorphism.

Note that since $X$ is smooth and finite type we know that $k[X]$ is normal.  Whether or not $k[X]^U$ is finitely generated, it is equal to the intersection of $k[X]$ and $k(X)^U$ in $k(X)$ and is hence by \cite{BComm} Example VII.1.3 (4) a Krull domain.  By \cite{BComm} VII.1.6 Theorem 4, $k[X]^U$ is therefore equal to the intersection of its localizations at height $1$ prime ideals.

We will show that pull-back $j^*: k[X]^U
\longrightarrow \Gamma(X/U,\O_{X/U})$ is a $k$-algebra isomorphism
if and only if $X/U$ has complement of codimension $\geq 2$ in
$\Spec k[X]^U$.  If $X/U$ has complement of codimension $\geq 2$ in
$\Spec k[X]^U$, any element $f \in \Gamma(X/U,\O_{X/U})$ extends uniquely to $k[X]^U$.  Indeed, this follows because any such $f$ can be
viewed as an element of every localization of $k[X]^U$ at a height
$1$ prime ideal. Conversely, assume $X/U$ has complement of
codimension $1$.   This means that the complement of $X/U$ in $\Spec
k[X]^U$ contains a non-trivial height $1$ prime ideal and hence
there exists an element $f \in \Gamma(X/U,\O_{X/U})$ not an element of
$k[X]^U$.

As in the proof of Theorem \ref{thm:affinequotient}, since the
quotient morphism is a faithfully flat and affine morphism, and $X$
is affine, we can identify $H^i(X/U,\O_{X/U})$ with
$H^i(U,\Gamma(X,\O_X))$ (again use the \u Cech resolution just prior
to the proof of  Theorem \ref{thm:affinequotient} in \S \ref{ss:affinequotient}).  By taking $i =
0$, this identification gives an isomorphism $\Gamma(X/U,\O_{X/U})
\cong k[X]^U$.   Combining this with the previous paragraph we see that $X/U$ always has complement of codimension
$\geq 2$ in $\Spec k[X]^U$.
\end{proof}

\section{Quotients of invariant hypersurfaces}
\label{s:hypersurfaces}
In this section, we focus on the study of quotients of
$\ga$-invariant hypersurfaces in affine spaces of the form ${\mathbb
A}(W)$ where $W$ carries a $k$-rational representation
of $\ga$; we will say that the corresponding $\ga$-action on ${\mathbb A}(W)$ is {\em linear} (and similarly for more general groups).  Henceforth, we abuse notation and write $W$ for the
affine space ${\mathbb A}(W)$.  In this section we give a geometric
characterization (see Theorem \ref{thm:geomchar1}) of when a $\ga$-invariant
hypersurface $X$ in $W$ admits the structure of a $\ga$-torsor over
an affine or strictly quasi-affine base $X/\ga$. Furthermore, we
characterize the corresponding form of the polynomial $f$ defining
$X$; all of the information is encoded in invariant and covariant
theory for $SL_2$ (see Theorem \ref{thm:algchar}). As a simple consequence
we will see that $k[X]^{\ga}$ is finitely generated (see Corollary
\ref{cor:hypersurfacefg}) in this situation; however, for $X$ of
higher codimension, we note this finite generation statement fails
(see Remark \ref{rem:highercodimfgfails}).  In the case of strictly quasi-affine $X/\ga$ we determine the boundary locus in $\Spec
k[X]^{\ga}$.

\subsubsection*{Constructing everywhere stable actions}
Let us first describe the method by which we shall produce varieties with everywhere stable actions of unipotent groups.  Suppose $U$ is a unipotent group and $G$ is any reductive group for which $U$ is a closed subgroup.  Note that any affine $U$-variety $X$ embeds as a closed subscheme of a finite dimensional $k$-rational $U$-representation $W$.  As the embedding $f: X \longrightarrow W$ is $U$-equivariant and closed, it follows from Corollary \ref{cor:morphismproperty} that the induced morphism
$$
G *_U f: G *_U X \longrightarrow G *_U W
$$
is a closed immersion as well.

It can be shown that every everywhere stable affine $U$-variety can be $U$-equivariantly embedded as a closed subscheme of the scheme $W^s$ for some finite dimensional $U$-representation $W$, but as we will not use this fact, we do not prove it here.  Instead, the following lemma is our main tool in construction of everywhere stable actions.

\begin{lem}
Suppose $W$ is a linear $U$-representation and $X$ is a quasi-affine $U$-stable subscheme of $W$ contained in $W^s$.  Then the $U$-action on $X$ is everywhere stable.
\end{lem}

\begin{proof}
By Corollary \ref{cor:morphismproperty}, the $U$-equivariant quasi-affine morphism $f: X \hookrightarrow W$ induces a quasi-affine morphism $G *_U f: G *_U X \hookrightarrow G *_U W$ (which is also locally closed).  By the theorem hypothesis together with definition \ref{defn:stable}, we see that $G *_U X \hookrightarrow (G *_U W)^s$ where $(G *_U W)^s$ is the set of stable points for the $G$-action on $G *_U W$.  We may then apply Proposition 1.18 of \cite{GIT} to conclude that the intersection $G *_U X \cap (G *_U W)^s$ is contained in the locus $(G *_U X)^s$.  We can identify the first intersection with $G *_U (W^s \cap X)$ by the definition of stability for $U$.  By assumption, $W^s \cap X = X$ and hence $X = X^s$ so that the $U$-action on $X$ is everywhere stable by Definition \ref{defn:everywherestable}.
\end{proof}

We will henceforth abuse terminology and refer to any $U$-stable quasi-affine subvariety $X$ of $W^s$ for a linear $U$-representation $W$ as an {\em everywhere stable subvariety of $W$}.  In the sequel, we will be interested in {\em everywhere stable hypersurfaces}, i.e., $U$-stable codimension $1$ closed subvarieties of $W$ contained in $W^s$.

\subsubsection*{Linear Representations}
Consider $\ga$ as a closed subgroup of $SL_2$ via the homomorphism defined by sending $x \mapsto \begin{pmatrix} 1 & 0 \\ x & 1 \end{pmatrix}$.  Given a $k$-rational representation $W$ of $\ga$, linearity has the following extremely useful consequence.

\begin{thm}[Jacobsen-Morozov (see e.g. \cite{AnKa} Theorem 19.5.1)]
\label{thm:jacobsenmorozov}
Given a linear $\ga$-representation on a vector space $W$,
i.e. a morphism $\rho: \ga \longrightarrow GL(W)$, there exists a
factorization of $\rho$ as the composite $\ga \hookrightarrow SL_2
\stackrel{\rho'}{\longrightarrow} GL(W)$.  Briefly, we will say that any
linear $\ga$-representation $\rho$ extends to a linear representation of
$SL_2$.
\end{thm}

\begin{rem}
If the $\ga$-action on $W$ admits, in Jordan canonical form, a
single Jordan block, then this extension is canonical (once given
the closed embedding $\ga \hookrightarrow SL_2$ as above). If there
are multiple blocks, then there is a choice of relative scaling
amongst the blocks for the action of the maximal torus $\gm \subset
SL_2$ of diagonal matrices. It is important to note that the above
extension induces a bijection between the isomorphism classes of
indecomposable representations of $\ga$ and $SL_2$.  Fix a
$\ga$-representation $W$ and any pair of extensions to $SL_2$
representations.  The resulting $SL_2$-representations are
necessarily isomorphic, as they have the same underlying sets of
irreducibles, but this isomorphism is non-canonical.  For a much
more complete discussion of functoriality in the
Jacobsen-Morozov theorem, see \cite{AnKa} \S 19 and especially Thm.
19.5.1.
\end{rem}

As in \S \ref{s:unipotents}, consider the inclusion $\iota: W \hookrightarrow
SL_2 *_{\ga} W$. Since the $\ga$-action on $W$ extends to a
$SL_2$-action, the product of the projection $SL_2 *_{\ga} W
\longrightarrow SL_2/\ga$ and the action map $SL_2 *_{\ga} W
\longrightarrow W$ defines a canonical isomorphism $SL_2 *_{\ga} W
\cong SL_2/\ga \times W$. Furthermore, $SL_2/\ga$ is the complement
of $0$ in the affine space $V$ where $V$ is the standard
representation of $SL_2$.  The point
$e \in SL_2/\ga$ then corresponds to the vector $\begin{pmatrix}0 \\
1 \end{pmatrix}$ in $V$.

\begin{lem}
\label{lem:W-transfer}
The composite map
$$
W \hookrightarrow SL_2/\ga \times W \hookrightarrow V \times W.
$$
determines (by pull-back) an isomorphism $k[V \times W]^{SL_2} \isomto k[W]^{\ga}$.  In addition, $k[W]^{\ga}$ is a finitely generated $k$-algebra.
\end{lem}

\begin{proof}
The inclusion $SL_2/\ga \hookrightarrow V$ has complement of
codimension $2$, hence $SL_2/\ga \times W$ has complement of
codimension $2$ in $V \times W$.  By normality of $V \times W$, restriction determines an
isomorphism of $k$-algebras $k[V \times W] \isomto k[SL_2/\ga \times W]$, and hence gives an isomorphism after taking $SL_2$-invariants.  Finally, we get an
isomorphism $k[W]^{\ga} \cong k[SL_2/\ga \times W]^{SL_2}$ by  Borel transfer (see \S \ref{ss:transfer}).  To conclude we note that $k[V \times W]^{SL_2}$ is finitely generated by Nagata's theorem on finite generation of rings of invariants under reductive group actions.
\end{proof}

We now discuss stability in $W$.  By Proposition \ref{prop:stable-for-red}, stability for the $\ga$-action on $W$ is
independent of the factorization $\rho'$ of $\rho: \ga \longrightarrow GL(W)$ through $SL_2$.  Therefore, we may extend the $\ga$-representation on $W$ to an $SL_2$-representation in such a way that the maximal torus of diagonal matrices $\gm \subset SL_2$ acts with equal weights on each indecomposable summand of the representation $W$.

Note that the stable set for the $SL_2$-action on $V \times W$ can be determined with the Hilbert-Mumford numerical criterion (see \cite{GIT} Theorem 2.1).  In
fact, because the inclusion $SL_2/\ga \times W \hookrightarrow V
\times W$ is a quasi-affine morphism, it follows from Proposition
1.18 of \cite{GIT} that the geometric points of $(SL_2/\ga \times W) \cap
(V \times W)^s$ are a subset of the geometric points of $(SL_2/\ga \times
W)^s$.  This distinguishes a subset of stable points on $W$ via
restriction from the closed immersion $\iota: W \rightarrow SL_2/\ga \times
W$.  Again, this set of geometric points necessarily underlies an open subscheme of $W$ which we denote by $W^{\bar{s}}$.  This
definition, coupled with a Hilbert-Mumford criterion computation,
yields the following Lemma, which is sufficient for the uses of
stability of this paper.

\begin{lem}
\label{lem:stable-bar}
The set $W^{\bar{s}}$ is an open subscheme of $W^s$.  It is the
complement in $W$ of the union of linear coordinate subspaces
defined by the vanishing of any coordinate with positive weight
relative to the action of the central torus $\gm$ of $SL_2$.
\end{lem}

Suppose that $X$ is a closed $\ga$-stable subvariety of $W$.  Then
we obtain an induced morphism $SL_2 *_{\ga} X \hookrightarrow SL_2
*_{\ga} W$ such that the following diagram commutes:
$$
\xymatrix{
X \ar[r]\ar[d] & SL_2 *_{\ga} X \ar[d] & \\
W \ar[r] & SL_2 *_{\ga} W \ar[r]^{\sim}& SL_2/\ga \times W.
}
$$
In particular, we can view $SL_2 *_{\ga} X$ as a locally closed
subvariety of $V \times W$.

\begin{lem}
\label{lem:invarhyper} Suppose $G$ is any linear algebraic group with trivial character group (e.g. $G = \ga$ or $SL_2$).  Any $G$-invariant hypersurface $X \subset W$ is actually defined by a $G$-invariant polynomial.
\end{lem}

\begin{proof}
Let $G$ be any linear algebraic group whose character group is trivial.
We know that any height $1$ prime ideal in $k[W]$ is principal, so
choosing a generator $f$ gives a defining equation for $X$.  Now,
$G$ acts on $X$ and hence the ideal defining $X$ is actually
$G$-stable.  Furthermore, the ideal sheaf of a codimension $1$
scheme is actually an invertible sheaf which we denote $\O(f)$.  Since $Pic(W)$ is trivial, we note that $\O(f)$ is isomorphic to
$\O_W$.  Since the character group of $G$ is trivial, the line
bundle $\O(f)$ admits a unique $G$-equivariant structure.  This
means that the module of sections $f \cdot k[W]$ admits a unique
$G$-invariant structure.  It follows that the action of $G$ on
$f$ is trivial and hence $f$ is actually invariant.
\end{proof}

By the transfer isomorphism of Lemma \ref{lem:W-transfer}, any
$\ga$-invariant polynomial $f$ on $W$ determines an $SL_2$-invariant
polynomial $F$ on $V \times W$.  Henceforth, we will use lower case Roman letters to denote $\ga$-invariants in $W$ and capital Roman letters to denote the corresponding $SL_2$-invariants in $V \times W$.  We now discuss the relation between
the  geometry of the hypersurfaces defined by the vanishing of each
of these polynomials.

\begin{lem}
\label{lem:F=0}
The $\ga$-invariant $f$ is irreducible if and only
if the associated $SL_2$-invariant $F$ is irreducible. If $X$ is
irreducible and everywhere stable,\footnote{The stability of $X$ may
not be necessary here.} then the closure $\overline{SL_2 *_{\ga} X}
\subset V \times W$ is given by $F = 0$.
\end{lem}

\begin{proof}
Note that if the vanishing locus of $f$ is a reduced scheme, then the vanishing locus of $F$ has the same property.
Now, for the forward implication, suppose $f$ is irreducible and $F$ is
not, and let $F_1$ and $F_2$ be distinct irreducible factors of $F$.
Then the vanishing locus of each $F_i$ defines an $SL_2$-invariant hypersurface
$X_i$ in $V \times W$.  The hypersurface $X_i$ is smooth and so all Weil divisors are Cartier; because the character group of $SL_2$ is trivial, the $F_i$ must be $SL_2$-invariant, as per the argument in Lemma
\ref{lem:invarhyper}.  By the transfer isomorphism we then have $f = f_1 f_2$, which is
a contradiction.

Similarly, for the reverse direction, if $F$ is irreducible but we can write $f =
f_1 f_2$, for non-constant $f_i$, then by the proof of Lemma
\ref{lem:invarhyper} both $f_1$ and $f_2$ are $\ga$-invariants.  Again by
 the transfer isomorphism this means $F = F_1 F_2$, which is a contradiction.

If $X$ is irreducible and stable then $X/\ga$ is irreducible, since
$X/\ga$ is a categorical quotient and a categorical quotient of an
irreducible variety is irreducible.  The total space of the
$SL_2$-principal bundle $SL_2 *_{\ga} X \longrightarrow X / \ga$ is
irreducible as it contains a dense, Zariski open subset isomorphic
to $V \times SL_2$ where $V$ is an open subset of $X/\ga$ and hence
irreducible by irreducibility of $X/\ga$.

Then $\overline{SL_2 *_{\ga} X}$ is the closure of an irreducible
variety, and must be a codimension $0$ closed subvariety of $F = 0$,
i.e., a maximal dimensional irreducible component of $F = 0$. But by
the first part of this lemma, $F = 0$ is itself irreducible.
\end{proof}

Choose coordinates $u, v$ on $V$, and let $w_1, \ldots, w_n$ be
coordinates on $W$.  Then $F = \sum_{i,j} F_{i,j} u^i v^j$ where the
$F_{i,j}$ are regular functions of $\{w_1, \ldots, w_n \}$.  Note
that $f$ can be recovered by restricting $F$ to the subvariety
defined by setting $u  =1$ and $v = 0$.

\begin{defn}
\label{defn:boundary} We will say that $X, f$, or $F$ {\em misses
the boundary} if $F_{0,0}$ is a non-zero constant.  We will say $X,
f$, or $F$ {\em contains the boundary} if $F_{0,0} = 0$.  Otherwise
we will say $X, f$, or $F$ {\em intersects the boundary}.
\end{defn}

\begin{rem}
\label{rem:algclosure} Of course, this terminology is shorthand for
how the geometric points of $\overline{SL_2 *_{\ga} X}$ relate to
the geometric points of the ``boundary" $\setof{0} \times W$ in $V
\times W$.  We can use geometric statements about intersections with the boundary
over $\kbar$ together with Galois descent to deduce which one of the three cases in the
definition holds.  In the sequel we will do this without further
formal justification.
\end{rem}

By Lemma \ref{lem:F=0}, $\overline{SL_2 *_{\ga} X} \subset V \times
W$ is the hypersurface with defining equation $F = 0$.  Note that
because $F$ is $SL_2$-invariant and the $SL_2$ action restricts to
$\setof{0} \times W$, it follows that $F_{0,0}$ is $SL_2$-invariant.
Thus there is a natural way to decompose $f$, namely $f = F_{0,0} +
g$, a sum of an $SL_2$-invariant function ({\em a priori} possibly
zero) and a $\ga$-invariant function that contains the boundary.
Everywhere stability of $X$ constrains $f$ further.

\begin{lem}
Let $X$ be everywhere stable.  Then $f = F_{0,0} + g$, where
$F_{0,0}$ is an $SL_2$-invariant not vanishing at the origin (i.e.,
with a nonzero constant term), and $g$ is a $\ga$-invariant that
contains the boundary.  In particular $f$ does not contain the
boundary.
\end{lem}
\begin{proof}
By the preceding discussion, all that need be shown is that
$F_{0,0}$ has a non-zero constant term.  But if it did not, then the
origin in $W$, which is not stable for any representation (indeed,
all homogeneous invariants automatically vanish there), would be a
solution to $f = 0$ and hence would be a point of $X$.
\end{proof}

\begin{thm}[Geometric characterization I]
\label{thm:geomchar1}
Let $X$ be an everywhere stable hypersurface in a linear representation $W$.  Then the quotient $X/\ga$ is affine if and only if $X$ misses the boundary as per Definition \ref{defn:boundary}. Otherwise,
\begin{itemize}
\item[i)]$X$ intersects the boundary and,
\item[ii)]$X/\ga$ is strictly quasi-affine and can be realized as an open subset of an affine variety with complement of codimension $\geq 2$.
\end{itemize}
\end{thm}

\begin{proof}
Note that it suffices to prove the result under the following two
assumptions: 1) $X$ is irreducible and 2) $k$ is algebraically
closed.  For the second statement, we note that stability is
preserved under base-change to the algebraic closure (see \cite{GIT}
Proposition 1.14; also see Remark \ref{rem:algclosure}).

Now, if $F$ misses the boundary, then $SL_2 *_{\ga} X$ is a closed
affine subvariety of $V \times W$.  Because $SL_2$ is reductive, the
quotient $SL_2 *_{\ga} X / SL_2 = X/\ga$ is affine.

Conversely, assume $X/\ga$ is affine yet $F$ intersects the
boundary.  Since $X$ is everywhere stable $F_{0,0}$ is non-constant.
In particular, the intersection $F_{0,0} \cap \setof{0} \times W$
has codimension $2$ in $\setof{F = 0} = \overline{SL_2 *_{\ga} X}$.
Let
$$\pi: \overline{SL_2 *_{\ga} X} \longrightarrow \overline{SL_2 *_{\ga} X}/\!/ SL_2$$
denote the categorical quotient map.  Then, again because $SL_2$ is
reductive, the image of $\pi$ is an affine variety.  Let $B$ denote
the complement of $SL_2 *_{\ga} X$ in $\overline{SL_2 *_{\ga} X}$.
By upper-semi-continuity of
the dimension of the fibers of $\pi$ (see e.g. \cite{Borel} Corollary 10.3), the scheme theoretic image of
$B$ under $\pi$ must have codimension at least two in the quotient
$\overline{SL_2 *_{\ga} X}/\!/SL_2 $.  Furthermore, $\pi(B)$ is
disjoint from $\pi(SL_2 *_{\ga} X)$ because $SL_2 *_{\ga} X \subset
\overline{SL_2 *_{\ga} X}^s$.  Consequently, $SL_2 *_{\ga} X /SL_2 =
X/\ga$ is an open subset of $\overline{SL_2 *_{\ga} X}/\!/SL_2$ with
complement codimension at least two.  Hence $X/\ga$ is a strictly
quasi-affine subvariety of an affine variety with complement of
codimension $\geq 2$.
\end{proof}

\begin{rem}
We know by Corollary \ref{cor:main} that $X/\ga$ is affine if and
only if $H^1({\mathfrak g}_a, k[X]) = 0$.  Equivalently, by Example
\ref{ex:imderivation}, if $D$ denotes the locally nilpotent
derivation defining the ${\mathfrak g}_a$-action on $k[X]$, we know
that $1 \in Im(D)$.  In practice, to determine whether $1 \in
Im(D)$, one first needs a generating set for $k[X]^{\ga}$.  If we
know in advance that $X/\ga$ is affine, then in principle one can
compute such a generating set (see Example \ref{ex:imderivation} and
references therein).

Although all $\ga$-invariants of $W$ restrict to $\ga$-invariants of
$X$, not all $\ga$-invariants of $X$ need extend to invariants of
$W$.  So one cannot {\em a priori} inherit the generating set of
invariants from $W$.  Nevertheless the next result, Corollary
\ref{cor:hypersurfacefg}, shows that at least assuming a normality hypothesis, the invariants do extend.
\end{rem}

\begin{cor}
\label{cor:hypersurfacefg} Suppose $X$ is a normal, everywhere
stable, hypersurface in $W$.  Then the ring of invariants
$k[X]^{\ga}$ is finitely generated, and in fact, all elements of
$k[X]^{\ga}$ arise as restrictions of elements of $k[V \times W]$.
\end{cor}

\begin{proof}
First, we show that if $X$ is normal, then so is $SL_2 *_{\ga} X$.
To see this, note that the morphism $\pi: SL_2 *_{\ga} X
\longrightarrow SL_2/\ga$ is a Zariski locally trivial fiber bundle
with fibers isomorphic to $X$.  Now, normality is a local property
for the Zariski topology, so picking an affine open cover of
$\setof{V_i}_{i \in I}$ of $X$ over which $\pi$ trivializes, we are
reduced to showing that $V_i \times X$ is normal, but this follows
because each $V_i$ is smooth.

Now, consider the inclusion $SL_2 *_{\ga} X \hookrightarrow
\overline{SL_2 *_{\ga} X}$.  By Lemma \ref{lem:F=0}, $\overline{SL_2
*_{\ga} X}$ is defined as the vanishing locus of $F$ in $V \times
W$.  Denote by $\widehat{SL_2 *_{\ga} X}$ the normalization of
$\overline{SL_2 *_{\ga} X}$.  Since $\widehat{SL_2 *_{\ga} X}$ is
affine and normal, the GIT quotient $\psi: \widehat{SL_2 *_{\ga} X}
\longrightarrow \widehat{SL_2 *_{\ga} X}/\!/SL_2$ is affine and
normal.  We now have the commutative diagram
$$
\xy
\xymatrix"*"{%
SL_2 *_{\ga} X/\!/SL_2 \ar[r] \ar[dr]& \widehat{SL *_{\ga} X}/\!/SL_2\ar[d] \\
 & \overline{SL_2 *_{\ga} X} /\!/SL_2} %
\POS-(30,15)
\xymatrix{%
SL_2 *_{\ga} X \ar["*"]^{q} \ar[r]\ar[dr]& \widehat{SL_2*_{\ga} X} \ar["*"]^{\hat{q}} \ar[d]\\
 & \overline{SL_2 *_{\ga} X} \ar["*"]^{\bar{q}}}%
\endxy
$$

Because $X$ is everywhere stable, $X/\ga$ is isomorphic to an open subscheme of $\overline{SL_2 *_{\ga} X} /\!/ SL_2$ with complement of codimension at least $2$, as argued in the proof of Theorem \ref{thm:geomchar1}.
Furthermore $X/\ga$, being a geometric and hence categorical
quotient of a normal variety, is itself normal.  Now, normalization
is the identity on normal varieties hence $X/\ga$ is an open subset
of $\widehat{SL_2 *_{\ga} X}/\!/SL_2$.  Furthermore, normalization
is a finite map, so the inverse image of a codimension $\geq 2$
subvariety has again codimension $\geq 2$; consequently, $X/\ga$ has
complement codimension $\geq 2$ in $\widehat{SL_2 *_{\ga}
X}/\!/SL_2$ .  By normality, all functions on $X/\ga$ extend
uniquely to functions on $\widehat{SL_2 *_{\ga} X}/\!/SL_2$, which
being affine implies $k[X/\ga]$ is finitely generated.  By the
properties of a geometric (indeed, categorical) quotient, $k[X/\ga]
= k[X]^{\ga}$, so $k[X]^{\ga}$ is finitely generated.

Note that if $F = 0$ is already a normal variety, then all of the
$SL_2$-invariant functions on $SL_2 *_{\ga} X$ extend to
$SL_2$-invariant functions on $\overline{SL_2 *_{\ga} X}$ which then
extend, by reductivity, to $SL_2$-invariant functions on $V \times
W$.  Hence, in this case, $k[X]^{\ga} \cong k[W]^{\ga}/(I(X) \cap
k[W]^{\ga})$.
\end{proof}

\begin{rem}
\label{rem:highercodimfgfails} If $X$ were of higher codimension in
$W$ then the conclusion of Corollary \ref{cor:hypersurfacefg},
namely that $k[X]^{\ga}$ need be finitely generated, does not
necessarily hold.  We now produce an explicit counter-example in the situation where $X$ is of codimension $\geq 3$ in $W$.  Take an action of
$\ga$ on an affine variety $X$ with non-finitely generated ring of
invariants.  Pick a $\ga$-equivariant embedding $X \hookrightarrow
W$.  Consider $SL_2$ as a hypersurface, necessarily $\ga$-invariant,
given by the defining equation $det(x) = 1$ in $End(V)$ (where $V$
is the $2$-dimensional representation of $SL_2$).  Note that this
hypersurface is everywhere stable.  Then $SL_2 \times X$ viewed as a
closed subvariety of $End(V) \times W$ is everywhere stable as well
and of codimension at least $2$.

By transfer, $k[X]^{\ga} = k[SL_2 *_{\ga} X]^{SL_2}$.  If $k[SL_2
*_{\ga} X]$ were finitely generated then by Nagata's theorem $k[SL_2
*_{\ga} X]^{SL_2}$ would be finitely generated as well.  Since $SL_2
\times X$ is everywhere stable, $k[SL_2 \times X]^{\ga}$ is
necessarily isomorphic to $k[SL_2 *_{\ga} X]$, but this is not
finitely generated.

Consider the example due to Daigle-Freudenberg of a $\ga$ action on
$\mathbb{C}^5$ with a non-finitely generated ring of invariants (see \cite{DaFr}). It is easy to recover their example as the $\ga$ action on
the closed subvariety $X$ in $W = \Sym^3 (V) \oplus V \oplus k$,
with coordinates $w_1, \ldots, w_7$, defined by the equations $w_7^3
= w_1, w_7^2 = w_5$. Thus $X \times SL_2$ can be presented as a
codimension $3$ closed invariant subvariety in $\Sym^3(V) \oplus
V^{\oplus 3} \oplus I$, which by the preceding argument cannot have
a finitely generated ring of invariants.
\end{rem}

In fact, $X$ being everywhere stable imposes further constraints.
Geometrically, a ``translate" of $X$ must contain all of the
non-stable points of $W$.  Let $c(f)$ denote the constant term of
$f$ (relative to the natural grading on $k[W]$ preserved by the
linear action of $\ga$).

\begin{lem}
\label{lem:WfinWs}
Let $X$ be everywhere stable.  Let $f^{\prime} = f - c(f)$. Then $W_{f^{\prime}} \subset W^s$.
\end{lem}

\begin{proof}
In fact we prove the {\em a priori} stronger statement that
$W_{f^{\prime}} \subset W^{\bar{s}}$.  Assume that some $w \in W$ is
simultaneously not in $W^{\bar{s}}$ and does not satisfy $f^{\prime}
= 0$.  Let $Y$ denote an irreducible component of $W \setminus W^s$
containing $w$. Any such $Y$ is a linear subspace of $W$, so in
particular contains the origin.  But $f^{\prime}$ vanishes at the
origin.  Thus $f^{\prime}$ is non-constant on $Y$ because $f'(w)$ is
non-zero. Since $Y$ is a closed and hence affine subvariety of $W$,
$f^{\prime}$ takes every value in $k$, in particular $c(f)$.  Thus
$f(y) = 0$ for some $y \in Y$, which violates the definition of
everywhere stable.
\end{proof}

\begin{thm}[Geometric characterization II]
\label{thm:geomchar2} Suppose $X$ is an everywhere stable
hypersurface in a linear representation $W$ with defining polynomial
$f$.  Then the quotient $X/\ga$ is an affine variety if and only if
$W_{f^{\prime}}/\ga$ is affine, where $f^{\prime} = f - c(f)$ for
$c(f)$ the uniquely defined constant term of $f$, in which case $X
\subset W_{f^{\prime}}$.
\end{thm}

\begin{proof}
By assumption, $X$ is contained in $W^s$.  But we also know $X
\subset W_{f'}$ and $W_{f'}$ is contained in $W^s$ by Lemma
\ref{lem:WfinWs}.  Therefore, a principal bundle quotient $W_{f'}
\longrightarrow W_{f'}/\ga$ exists as a scheme.

Now, suppose $W_{f'}/\ga$ is affine.  Since $X \subset W_{f'}$ is a
$\ga$-invariant closed subscheme, it follows that $X/\ga$ is closed
in $W_{f'}/\ga$ (since the $X/\ga$ is a geometric quotient).
Therefore, $X/\ga$ is affine.

Conversely, assume that $X/\ga$ is an affine variety.  Since $X$ is
everywhere stable, the defining polynomial $f$ must have non-zero
constant term (the origin is unstable).  If $f' = f - c(f)$, by
Lemma \ref{lem:WfinWs} we know that $W_{f'} \subset W^s$. Therefore,
it suffices to show that $f$ does not intersect the boundary.  If $f$ misses the boundary, then $f'$ contains the boundary and so
$W_{f'}$ must miss the boundary.  Therefore, $W_{f'}/\ga$ must be
affine by Theorem \ref{thm:geomchar1}.
\end{proof}

We now give a somewhat ``geometric" interpretation for a function
that contains the boundary.

\begin{lem}\label{lem:geom-contain-boundary}
Let $g$ be a $\ga$-invariant function on $W$.  If $g$ contains the
boundary then $W_g \subset W^s$.
\end{lem}

\begin{proof}
Again we prove the {\em a priori} stronger statement that $W_g
\subset W^{\bar{s}}$.  By Lemma \ref{lem:W-transfer}, $g$ corresponds to an
$SL_2$-invariant function $G$ on $V \times W$.  Recall that we used
the notation $u,v$ for the coordinate functions on $V$, with $u$
having positive weight and $v$ negative weight for the action of the
torus $\gm$ contained in the normalizer of $\ga$ in $SL_2$. Because
$g$ contains the boundary, each term of $G$ must contain a factor of
$u$ or $v$. Since $G$ is an $SL_2$-invariant it is in particular a
$\gm$-invariant, so each term of $G$ must be weight $0$. However, $g
= G |_{u=0, v=1}$, and hence must consist of terms that are of
strictly positive weight for the $\gm$ action on $W$.  By Lemma
\ref{lem:stable-bar} the non-stable set is contained in the locus of
points where all coordinate functions of strictly positive weight,
and hence all the terms of $g$, vanish.
\end{proof}

For open affines of the form $W_h$ the cohomological vanishing
criterion of Theorem \ref{thm:affinequotient} has a very simple
interpretation.

\begin{lem} \label{lem:h^k_im(D)}
Let $h$ be a $\ga$-invariant function on $W$.  The quotient
$W_h/\ga$ is affine and $W_h \longrightarrow W_h/\ga$ is a trivial
principal bundle if and only if some power of $h$ lies in $Im(D)$.
\end{lem}

\begin{proof}
By assumption we have $h \in Ker(D)$ since $h$ is a $\ga$ invariant.
By Theorem \ref{thm:affinequotient}, the quotient map $W_h
\longrightarrow W_h/\ga$ is a trivial principal bundle over the
affine variety $W_h/\ga$ if and only if $H^1({\mathfrak g}_a,k[W_h])
= 0$.  Example \ref{ex:imderivation} shows that this happens if and
only if there exists a function $s \in k[W_h]$ such that $D(s) = 1
\in k[W_h]$.

Let us first show that if some power of $h$ lies in $Im(D)$, that we
can find an $s$ such that $D(s) = 1$.  Indeed, suppose $h^k = D(g)$.
If we set $s = g/h^k$; then, by the quotient rule, $D(s) = 1$.

Conversely, assume there exists an $s \in k[W_h] = k[W][h^{-1}]$ such that $D(s) = 1$.  In coordinates,
$$s = g_0 +
\frac{g_1}{h} + \frac{g_2}{h^2} + \ldots \frac{g_k}{h^k} = \frac{h^k
g_0 + h^{k-1} g_1 + \cdots h g_{k-1} + g_k}{h^k}.
$$
Then, again by the quotient rule, $D(s) = D(h^k g_0 + h^{k-1} g_1 +
\cdots h g_{k-1} + g_k)/h^k$.  If we have $D(s) = 1$, then $D(h^k
g_0 + h^{k-1} g_1 + \cdots h g_{k-1} + g_k) = h^k$; in other words,
$h^k \in Im(D)$.
\end{proof}

\begin{lem}
A $\ga$-invariant function $f^{\prime}$ contains the boundary if and
only if there exists a positive integer $k$ such that
$(f^{\prime})^k \in Im(D) \cap Ker(D)$.
\end{lem}

\begin{proof}
Since $f^{\prime}$ contains the boundary,
$$SL_2 *_{\ga} W_{f^{\prime}} = ((V \setminus \setof{0}) \times W)_{F'} = (V \times W)_{F^{\prime}}.$$
In particular $SL_2 *_{\ga} W_{f^{\prime}} $ is affine.  The GIT
quotient of an affine variety by a reductive group is affine, so
$SL_2 *_{\ga} W_{f^{\prime}}/SL_2 = W_{f^{\prime}}/\ga$ is therefore
affine.  By definition of stability and Lemma \ref{lem:h^k_im(D)},
it follows $(f^{\prime})^k \in Im(D) \cap Ker(D)$.

Let $(f^{\prime})^k = D(g)$. As in the first argument in Lemma
\ref{lem:h^k_im(D)}, let $s = g/(f^{\prime})^k$, so that $D(s) = 1$.
Thus $W_{(f^{\prime})^k} = W_{f^{\prime}} \longrightarrow
W_{f^{\prime}}/\ga$ is a trivial principal $\ga$-bundle over the
affine variety $W_{f'}/\ga$, and $W_{f^{\prime}} \subset W^s$.

Now, we have the equality $W_{f^{\prime}}/\ga \cong ((V \setminus
\setof{0}) \times W)_{F^{\prime}}/SL_2$.  But, because the inclusion
$(V \setminus \setof{0}) \times W \subset V \times W$ is open with
complement of codimension 2 and since $V \times W$ is normal, all
regular functions extend.  Since the inclusion is
$SL_2$-equivariant, this extension sends $SL_2$-invariant regular
functions to $SL_2$-invariant regular functions.  Then $(V \times
W)_{F^{\prime}}//SL_2$ is an affine scheme and its coordinate ring
is identified with the coordinate ring of $((V \setminus \setof{0})
\times W)_{F^{\prime}}/SL_2$.  Therefore, we can identify the
quotients $(V \times W)_{F^{\prime}}//SL_2$ and $((V \setminus
\setof{0}) \times W)_{F^{\prime}}/SL_2$.

It follows that any point $x \in (\{ 0 \} \times W)_{F^{\prime}}$
has to be in the closure of an orbit $\mathcal{O}$ of $((V \setminus
\setof{0}) \times W)_{F^{\prime}}$, which violates the definition of
stability of $\mathcal{O}$.  This implies $(\setof{0} \times
W)_{F^{\prime}}$ is in fact empty and therefore that $f^{\prime}$
contains the boundary.
\end{proof}

Combining these lemmas, we are led to a simple algebraic characterization of everywhere stability.

\begin{thm}[Algebraic characterization]
\label{thm:algchar}
Suppose $X$ is a $\ga$-invariant hypersurface defined by a $\ga$-invariant polynomial $f$.  Then $X$ is everywhere stable if and only if there is a decomposition $f = F_{0,0} + g$ satisfying
\begin{itemize}
\item[i)] Either the function $F_{0,0}$ is a non-zero constant or $F_{0,0}$ is a {\em stable} $SL_2$-invariant with non-zero constant term (we will say that $F_{0,0}$ is a stable $SL_2$-invariant if $F_{0,0}(x)$ is non-vanishing only if $x \in W^s$), and
\item[ii)] the function $g$ is a $\ga$-invariant that contains the boundary, and some positive power of $g$ lies in $Im(D)$.
\end{itemize}
Moreover, the decomposition $f = F_{0,0} + g$ is unique.  Furthermore, in the above situation, the quotient $X/\ga$ is affine if and only if $F_{0,0}$ is constant.
\end{thm}

\begin{rem}
In fact, these lemmas show that {\em any} $\ga$-invariant function
admits a natural decomposition into an $SL_2$ invariant function
$F_{0,0}$ and a $g$ such that $g^k \in Im(D) \cap Ker(D)$ for some
positive integer $k$.  The everywhere stable condition simply
constrains $F_{0,0}$ to be a ``stable" invariant.  It would be
interesting to know whether this decomposition of $\ga$ invariants
can be deduced from classical $SL_2$-covariant theory.
\end{rem}

When enough ``ingredients" in this picture are smooth, then the
quotients have a particularly clean presentation as open subsets of
hypersurfaces in $\Spec k[W]^{\ga}$.

\begin{lem} \label{lem:F00-smooth}
Assume that the $\ga$-invariant hypersurface $X$ in $W$, defined by
$f = 0$, is smooth. If the variety in $W$ defined by $F_{0,0} = 0$
is smooth, then the variety in $V \times W$ defined by $F = 0$ is
smooth.
\end{lem}

\begin{proof}
We can assume that $k$ is algebraically closed as $X$ is smooth if and only if its base extension to $\kbar$ is non-singular.  If $X$ is a smooth affine scheme, it follows by Corollary
\ref{cor:morphismproperty} and the fact that $SL_2/\ga$ is smooth that the induced scheme $SL_2 *_{\ga} X$
is smooth as well.  Let us fix coordinates $u,v$ on $V$ as we have
above.  Now, the singular locus of the hypersurface defined by the
polynomial $F$ must lie in the locus of points where $u = 0$ and $v
= 0$, i.e. in $\setof{0} \times W$.  We may now apply the Jacobian
criterion for smoothness.

The singular locus of the hypersurface defined by $F$ is the
simultaneous vanishing locus of the partial derivatives
$\frac{\partial F}{\partial v}$, $\frac{\partial F}{\partial u}$,
and $\frac{\partial F}{\partial w_i}$ for all $i$.  Now,
$\frac{\partial F}{\partial w_i} = \frac{\partial F_{0,0}}{\partial
w_i} + u( \cdots) + v(\cdots)$.  It follows that any singular points
must be singularities of $\frac{\partial F_{0,0}}{\partial w_i} =
0$.
\end{proof}

\begin{thm} \label{thm:F00-smooth}
Let $X$ be a smooth everywhere stable hypersurface in $W$ defined
by the vanishing of a $\ga$-invariant polynomial $f$.  Let $F_{0,0}$ be the $SL_2$-invariant polynomial appearing in the decomposition of $f$ as in Theorem \ref{thm:algchar}.  Assume the variety $Y$ in $W$ defined by $F_{0,0}$ is smooth.  Then the quotient $X/\ga$ is an open subvariety of the normal variety $\Spec k[W]^{\ga}/(f)$.  Furthermore, the boundary of $X/\ga$ in $\Spec k[W]^{\ga}/(f)$ is of codimension $2$ and is identified with the quotient $Y/SL_2$, which has at most finite quotient singularities.
\end{thm}

\begin{proof}
By Lemma \ref{lem:F00-smooth} the variety defined by the vanishing of $F$, which we denote by $\overline{SL *_{\ga} X}$, is smooth.  Hence any $SL_2$-invariant function on $SL_2 *_{\ga} X$ extends to an $SL_2$-invariant function on $\overline{SL_2 *_{\ga} X}$.  Any $SL_2$-invariant function on an affine subvariety extends to an $SL_2$-invariant function on the ambient affine space; here this is given by $V \times W$.

Consequently, we have isomorphisms $k[X]^{\ga} \cong k[V \times W]^{SL_2}/(F) \cong k[W]^{\ga}/(f)$.  Thus the quotient $X/\ga = SL_2 *_{\ga} X/SL_2$ is an open subscheme of $\overline{SL_2 *_{\ga} X}/\!/SL_2 = \Spec (k[V \times W]/(F))^{SL_2} \cong \Spec k[V \times
W]^{SL_2}/(F) \cong \Spec k[W]^{\ga}/(f)$.  This latter variety is normal because $\overline{SL *_{\ga} X}$ is smooth, hence normal, so the categorical quotient $\overline{SL_2 *_{\ga} X}/\!/SL_2$ is normal.

By Theorem \ref{thm:algchar}, $F_{0,0}$ is a stable $SL_2$-invariant. Thus the boundary $F_{0,0} = 0$, and hence $F = 0$ is in fact everywhere stable.  Thus, we have an identification $\overline{SL_2 *_{\ga} X}/\!/SL_2 = \overline{SL_2 *_{\ga} X}/SL_2$, so in particular the boundary of the quotient is isomorphic to the quotient $F_{0,0}/SL_2$.
\end{proof}

\begin{rem}
\label{rem:Stein}
For this remark, let us assume that $X$ is defined over $\cplx$ and let $X(\cplx)$ denote the associated analytic space.  It follows from Theorem \ref{thm:F00-smooth} that when the quotient $X/\ga$ is quasi-affine, the analytic space $X/\ga(\cplx)$ is not Stein because it is the complement of a codimension $2$ analytic subspace in a normal Stein space.
\end{rem}

It is thus easy to produce examples of affine or strictly
quasi-affine quotients; indeed, any possible example can be produced
through the $SL_2$-invariant (and covariant) theory of $W$.

\subsubsection*{A quasi-affine quotient}
Let us illustrate the algebraic characterization by giving the
simplest example of a strictly quasi-affine quotient of an
everywhere stable action.

\begin{ex}
Consider $SL_2$ as a hypersurface $X$ in $W = V \oplus V$ determined
by the vanishing of the polynomial $f = 1 - w_0w_3 + w_1w_2$.  Note
that $f$ is a stable-$SL_2$-invariant with non-zero constant term
because the unstable locus is the set of all points $w_0 = w_2 = 0$.
Therefore, by Theorem \ref{thm:algchar}, the hypersurface defined by
$f$ is everywhere stable and has a quasi-affine quotient.  Here, the
invariants are generated by $w_0$ and $w_2$ and the image $X/\ga$ is
isomorphic to ${\mathbb A}^2 \setminus \setof{0}$.

Similarly, if $\phi$ is a degree $d$ function of one variable with
no repeated roots and no constant term, then let $f = 1 -
\varphi(w_0w_3-w_1w_2)$.  The boundary in $F = 0$ consists of $d$
disjoint isomorphic copies of $SL_2$.  The quotient is identified
with $\mathbb{A}^2 \setminus \setof{p_1, \ldots, p_d}$, where $p_i$
are points representing the boundary components.  Since the action
of the automorphism group $Aut_k({\mathbb A}^2)$ on $\mathbb{A}^2$
is $d$-transitive for any positive integer $d$, it follows that, up
to isomorphism, the complement in ${\mathbb A}^2$ of any set of
$d$-points may always be realized as a $\ga$-quotient.
\end{ex}

\section{Examples: $\aone$-contractible strictly quasi-affine varieties}
\label{s:examples}
In this section, we use Theorem \ref{thm:algchar} to give many
examples of $\aone$-contractible smooth varieties in every dimension
$\geq 4$.  In order to do this, we will construct everywhere stable
actions of $\ga$ on affine spaces with strictly quasi-affine
quotients.

\subsubsection*{$SL_2$-representations}
Let $V$ denote the standard $2$-dimensional representation of
$SL_2$. Since every linear representation of $\ga$ extends to an
$SL_2$-representation, we will index representations of $\ga$ by the
corresponding representations of $SL_2$.  Fix coordinates
$w_0,\ldots,w_k$ on $\Sym^k V$ corresponding to a basis of weights
with $w_0$ being $\ga$-fixed.  Let $\partial_i$ denote the
derivation $\frac{\partial}{\partial w_i}$.  In these coordinates,
the locally nilpotent derivation defining the $\ga$-action on
$\Sym^k V$ is given by
$$
D_{\Sym^k V} = \sum_{i=0}^{k-1} (k-i)w_{i} \partial_{i+1}.
$$
A general linear representation $W$ of $\ga$, extends to an $SL_2$
representation that decomposes as $W = \bigoplus \Sym^{k_j} V$.  The
locally nilpotent derivation $D$ determining the $\ga$-action on $W$
then takes the form  $$D = \sum_{i \in I} c_i w_i \partial_{i+1}$$
where the set $I$ and the coefficients $c_i$ are determined by the
$k_i$.

\subsubsection*{Everywhere stable embeddings of affine space}
We are interested in realizations of ${\mathbb A}^n$ as an
everywhere stable hypersurface in a linear representation $W$. In
order to do this, choose coordinates $z_1,\ldots,z_n$ on ${\mathbb
A}^n$, and $w_0,\ldots,w_n$ on $W$.  As usual, an embedding
${\mathbb A}^n \hookrightarrow W$ is determined by specifying
polynomials $w_i = w_i(z_1,\ldots,z_n)$.  The embedding being
$\ga$-equivariant is equivalent to the locally nilpotent derivation
$D$ defining the $\ga$-action on $W$ restricting to a locally
nilpotent derivation on ${\mathbb A}^n$; by abuse of notation, we
will also call this restricted derivation $D$.

Let $f$ be the polynomial relation between $w_i$; this defines a
hypersurface $X$ in $W$.  By Lemma \ref{lem:invarhyper}, the
polynomial $f$ is necessarily $\ga$-invariant and hence $D(f)$ must
vanish identically.  This imposes the condition
$$
\sum_{i \in I} c_i w_i \partial_{i+1}f = 0,
$$
and $I$ is an index set as above.

We now construct a class of $\ga$-equivariant embeddings of
${\mathbb A}^n$ into $W$ as follows.  Consider the morphism
${\mathbb A}^n \longrightarrow W$ defined by setting $w_i = z_i$ for
all $i \neq 0$, and let $w_0 =h(z_1,\ldots,z_n)$ where $h$ is a
$\ga$-invariant function; observe that this morphism is actually an
embedding.  The corresponding hypersurface $X$ in $W$ is then the
vanishing locus of the polynomial $f = h(w_1,\ldots,w_n) - w_0$.  By
Theorem \ref{thm:algchar}, in order that $X$ be an everywhere stable
hypersurface with quasi-affine quotient, the invariant $h$ must
decompose as $F_{0,0} + g$ where i) $F_{0,0}$ is a stable
$SL_2$-invariant with non-vanishing constant term and ii) some
strictly positive power of $g$ is in the image of $D$.

\begin{thm}
\label{thm:dimension4}
For every $m \geq 4$, there exists a
denumerably infinite collection of pairwise non-isomorphic $m$-dimensional exotic $\aone$-contractible varieties, each admitting an embedding into a
smooth affine variety with pure codimension $2$ smooth boundary.
\end{thm}

\begin{proof}
Consider the representation $W = \Sym^q \oplus \Sym^{2p+1}(V) \oplus W'$, where $p, q \geq 1$ and $W'$ is any linear $\ga$
representation.  Choose coordinates $w_0,\ldots,w_{2p+q+ \dim W'
+2}$ which are a basis of $\gm$-eigenvectors for an $SL_2$-action on
$W$ extending the given $\ga$-action; by convention, the first
$q+1$-coordinates will be coordinates on $\Sym^{q}(V)$, the next
$2p + 2$-coordinates will be coordinates on $\Sym^{2q+1}(V)$, and the
remaining $w_i$ will be coordinates on $W'$.

Because $\Sym^{2p+1}(V)$ has an $SL_2$-stable point, it has an
$SL_2$ GIT quotient of dimension $2p-1 \geq 1$, so in particular
there exists a non-constant homogeneous $SL_2$ invariant $\Delta$ on
$W$ (for example, the discriminant). Furthermore, because
$\Sym^{2p+1}(V)$ is even dimensional, all of its points are either
stable or unstable (in that all homogeneous invariants vanish on the
unstable set); consequently all homogeneous $SL_2$ invariants are in
fact stable $SL_2$ invariants.  Let $w_{0}$ denote the invariant
coordinate in direct summand $\Sym^q(V)$ of $W$; note that $w_{0}
= D(w_{1})$ (and $w_{1}$ is an element of $\Sym^q(V)$ because $q \geq 1$).  In particular, this explains why we must choose $q \geq 1$.

Suppose $\varphi$ is a polynomial in one variable of strictly
positive degree, with no multiple roots, and satifying $\varphi(0) \neq -1$.  Let $h_\varphi = 1 + \varphi(\Delta)$. Note that $h_\varphi$ is a stable $SL_2$-invariant with non-vanishing constant term. Therefore, if $f = h_\varphi - w_{0}$, the discussion just prior to the theorem says that $f$ defines a hypersurface isomorphic to ${\mathbb A}^{m+1}$ in $W$, where $m = 2p+q+\dim(W')+1$. By Theorem \ref{thm:algchar} the quotient
$X_\varphi = {\mathbb A}^{m+1}/\ga$ is strictly quasi-affine.

We now claim that if $\varphi$ and $\varphi'$ are two one-variable
polynomials whose degrees differ, then the resulting quotients are
non-isomorphic.  Note that $X_\varphi$ is an open subvariety of the
affine variety $\Spec k[z_1,\ldots,z_n]^{\ga}$ (which is finitely
generated by Corollary \ref{cor:hypersurfacefg}) with complement of
codimension at least $2$.  Let us denote the boundary by
$B_\varphi$.

Suppose $\varphi$ and $\varphi'$ are a pair of strictly positive
degree polynomials in one variable with no repeated roots.
Let $d = deg(\varphi)$.  If $F_\varphi$ denotes the $SL_2$-invariant
hypersurface equation attached to $h_\varphi$, then the boundary of
the hypersurface $F_\varphi$ in $V \times W$ is given by
$F_{\varphi,0,0} = 0$.  Here $F_{\varphi,0,0} = \varphi(\Delta)$.
Hence the boundary $B_\varphi$ is the $SL_2$-quotient of the
vanishing locus of $F_{\varphi,0,0}$. This equips $B_\varphi$ with
the structure of vector bundle over the vanishing locus of
$h_\varphi$. Since $\varphi$ has $d$ distinct roots, the components
are of the form $\Delta = c_i$ where
$i$ runs from $1$ to $d$; these components are obviously disjoint.
Taking the $SL_2$-quotient, we see that $B_\varphi$ therefore has
$d$ distinct components as well.

Finally, note that because $W$ has no points of finite isotropy in
$SL_2$, the variety $(1+\Delta = 0)/SL_2$ is smooth, and for generic
$\varphi$ the variety $h_{\varphi}/SL_2$ is smooth.  Likewise, $V
\times W$ has no points of finite isotropy, so for generic $\varphi$
the variety $(F_{\varphi} = 0)/SL_2$ is smooth.  The smoothness
assertions in the theorem follow.
\end{proof}

\begin{rem}\label{rem:Winkelex}
Note that $W'$ plays no essential role in the above proof; it is simply present to allow for the construction of more examples.  Indeed, to produce yet more examples, the factor of $\Sym^{2p+1} V$ in $\Sym^{q} V \oplus \Sym^{2p+1} V \oplus W'$ can be replaced by any $\ga$-representation with at least one stable $SL_2$-invariant -- for example, $V \oplus V$ with its unique quadratic $SL_2$-invariant.  If one takes $q=1$ and $W'= \setof{0}$, then the hypersurface defined by the vanishing of $w_0 = 1 + (w_2w_5-w_3w_4)$ is a $\ga$-equivariant linear embedding of Winkelmann's example (see \cite{Wi} \S 2) and the quotient of course then agrees with his.
Specifically, by Corollary \ref{cor:hypersurfacefg} the
$\ga$-invariants for $X$ here are just restrictions from
$k[W]^{\ga}$, namely:
$$w_0, w_2, w_4, w_0 w_3 -w_1 w_2, w_0 w_5 -
w_1 w_4, w_2 w_5 - w_3 w_4.$$
Imposing the hypersurface equation
for $X$ there are only $5$ generating invariants, and one relation.
The relation gives the hypersurface equation $x_1 x_4 - x_2 x_3 -
x_5(x_5+1)=0$ in $\mathbb{A}^5$, where $x_1, \ldots, x_5$ are
identified, in order, with the last five invariants above.
\end{rem}

\begin{thm}\label{thm:moduli}
For every $m \geq 6$ and every $n
> 0$:
\begin{itemize}
\item there exists a connected $n$-dimensional scheme $S$ and a smooth
morphism $f: X \longrightarrow S$ of relative dimension $m$, whose
fibers over $k$-points are $\aone$-contractible and quasi-affine,
not affine, and pairwise non-isomorphic.

\item The morphism $f: X \longrightarrow S$ admits a partial compactification to a flat family $\bar{f}:
\overline{X} \longrightarrow S$ whose fibers over $k$-points are
smooth affine varieties. Furthermore, for any $k$-point $t \in S$,
the map $X_t \rightarrow \bar{X}_t$ is an open immersion with a
smooth complement of codimension $\geq 2$
\end{itemize}
\end{thm}

\begin{proof}
Let $W$ be a linear $\ga$-representation.  Suppose we choose a
family $f_t$ of polynomials in $W^s$ parameterized by some base
scheme $T$, such that the corresponding hypersurfaces $X_t$ are all
everywhere stable.  Let $F_{t}$ is the corresponding induced family
of polynomials in $G *_U W$ and let $F_{t,0,0}$ be the polynomial
defining the boundary.  Let $\overline{SL_2 *_{\ga} X_t}$ be the
vanishing locus of $F_t$ in $V \times W$.  Assuming $X_t$ are all
smooth, we know the quotients $X_t/\ga$ are all smooth. Furthermore,
by Lemma \ref{lem:F00-smooth}, if the vanishing loci of $F_{t,0,0}$
define  smooth subvarieties of $W$, it follows that $F_t$ are all
smooth varieties as well.  Denote these vanishing loci by $B_t$.

If the $F_t$ are smooth, then the categorical quotients
$\overline{SL_2 *_{\ga} X_t}/\!/SL_2$ are normal varieties.
Furthermore, if $B_t$ is smooth, then $B_t/\!/SL_2$ is a normal
variety as well.  Now, any isomorphism $\psi: X_t/\ga \isomto
X_{t'}/\ga$ extends to an isomorphism $\bar{\psi}: \overline{SL_2
*_{\ga} X_{t}} \isomto \overline{SL_2 *_{\ga} X_{t'}}$ by normality.  The
isomorphism $\bar{\psi}$ then would restrict to an isomorphism
$\hat{\psi}: B_t \isomto B_{t'}$.  Therefore, if $B_{t}/\!/SL_2$ and
$B_{t'}/\!/SL_2$ are non-isomorphic, the quotients $X_{t'}/\ga$ and
$X_{t'}/\ga$ are not isomorphic.

If $w_0,\ldots,w_n$ are coordinates on a linear representation $W$,
as in Theorem \ref{thm:dimension4}, we will use invariants of the
form $f_t = h(w_1,\ldots,w_n) - w_0$, where $h$ is a stable
$SL_2$-invariant.  Indeed, suppose we can find $\setof{\Delta_i}$,
$i = 1,\ldots,j$ a collection of stable $SL_2$-invariants in $W$
depending only on $w_1,\ldots,w_n$ that are algebraically
independent.   Then if we fix a family of polynomials $\varphi_t$ in
$j$ variables, such that $\varphi_t(0) \neq -1$, the function $h_t = 1 +
\varphi_t(\Delta_1,\ldots,\Delta_j) - w_0$ gives a family of
embeddings of affine space into $W$ parameterized by $t$.
Furthermore, as $\Delta_i$ are invariants and algebraically
independent, we can treat the $\Delta_i$ as coordinates on the
quotient $B_t/\!/SL_2$.  In other words, we can view $B_t /\!/SL_2$
as a closed subvariety of $W/\!/SL_2$ defined by the vanishing of
$\varphi_t(\Delta_1,\ldots,\Delta_j)+1$.

Assuming these families admit $n$-dimensional moduli, which  we
defer along with a discussion of smoothness to Lemmas
\ref{lem:families} and \ref{lem:smoothfamilies}, let us
complete the proof of the theorem. We can assume that $t$ takes
values in a non-singular curve $T$. Then, the assignment $X_t
\longrightarrow t$ defines an algebraic family of varieties $\pi: X
\longrightarrow T$ over a curve. Furthermore, the closures
$\overline{SL_2 *_{\ga} X_t}/\!/SL_2$ define an algebraic family of
normal schemes (see \cite{Ha}  Ch. III Defn. 9.10) in which $X_t$
are open of codimension $\geq 2$, in other words we get a morphism
$\bar{\pi}: \overline{SL_2 *_{\ga} X}/\!/SL_2 \longrightarrow T$
factoring $\pi$.  Such families are necessarily flat by \cite{Ha}
Ch. III Thm. 9.11. The family $X_t$ is in fact smooth as the
geometric fibers $X_t$ are all smooth. Therefore, we get a family
$X' \longrightarrow C$ which is smooth. Taking products of such
families, we can increase the dimension of the parameter space of
the family as we wish.
\end{proof}

\begin{rem}
This discussion recovers Winkelmann's example of a family of
quasi-affine quotients of $\mathbb{A}^7$ (see \cite{Wi} \S 4). There
is, however, a gap in Winkelmann's argument for existence of
families of quasi-affine quotients.  If $X$ and $X'$ are two
quasi-affine schemes which are realized as closed subschemes of
affine schemes $\overline{X}$ and $\overline{X'}$, then an
isomorphism between $X$ and $X'$ extends uniquely to the
normalizations of $\overline{X}$ and $\overline{X'}$.  Winkelmann
then claims that this extended isomorphism induces an isomorphism of
the normalizations of the boundaries $\overline{X'} \setminus X'$
and $\overline{X} \setminus X$.  However, this is not necessarily
true.  Indeed, it is possible that $\overline{X}$ is normal (so that
its normalization is trivial) and $\overline{X} \setminus X$ is not
normal.  Furthermore, even if $\overline{X} \setminus X$ were normal, unless
$\overline{X}$ is normal there may well be moduli of $\overline{X}
\setminus X$ all of which have the same inverse image in the
normalization of $\overline{X}$.  For us, all of these problems are
circumvented by application of Lemma \ref{lem:F00-smooth}.
\end{rem}

The remaining step in the proof of Theorem \ref{thm:moduli} is to
establish the intuitively ``clear" statement that there exist
boundaries $B_t/SL_2$, expressed as $\varphi_t(\Delta_1, \ldots,
\Delta_j) - 1=0$ in $\Spec k[W]^{SL_2}$, which admit many
deformations.  Rather than argue in complete generality, we now consider a simple special case sufficient for the theorem; the more general
formulation would follow exactly the same lines.

\begin{lem}\label{lem:families}
Let $W = V^{\oplus k}$, for $V$ the standard two-dimensional
representation of $SL_2$, and for $k \geq 4$.  Then for any $p$
there exists a $p$-dimensional family $X \longrightarrow T$ of
non-isomorphic smooth hypersurfaces in $\Spec k[W]^{SL_2}$ defined
by $\varphi_t(\Delta_1, \ldots, \Delta_j)+1=0$.
\end{lem}

\begin{proof}
Let $x_{i, \varepsilon}$, for $i \in\{ 1, \ldots, n \}$ and
$\varepsilon$ either $0$ or $1$, be coordinates on $W$; identify
them with the usual coordinates via $w_{2(i-1)+ \varepsilon} = x_{i,
\varepsilon}$.  The ring $k[W]^{SL_2}$ is generated by
$(\frac{k}{2})$ quadratic stable invariants: $x_{i,0}x_{j,1} -
x_{i,1}x_{j,0}$.  (That all $SL_2$ invariants in this representation
are stable follows from the Hilbert-Mumford criterion.)  Consider
the three invariants, call them $\Delta_1, \Delta_2, \Delta_3$,
associated with the last $3$ direct summands in $W = V^{\oplus k}$.
They are the coordinate functions of $V^{\oplus 3}//SL_2 \cong
\mathbb{A}^3$.

Consider a hypersurface $B_t$ in $W$ defined by $\varphi_t(\Delta_1,
\Delta_2, \Delta_3) + 1$; then $B_t \cong V^{\oplus k-3} \times
Y_t$, where $Y_t$ denotes the variety defined by
$\varphi_t(\Delta_1, \Delta_2, \Delta_3) + 1 = 0$ in $V^{\oplus 3}$.
Observe that $SL_2$ acts freely on any $SL_2$-invariant everywhere
stable hypersurface.  Note that the projection morphism $B_t
\longrightarrow Y_t$ equips $B_t$ with the structure of a trivial
$SL_2$-equivariant vector bundle over $Y_t$.   Since the
$SL_2$-action is free on $Y_t$, the projection morphism $B_t
\longrightarrow Y_t$ descends to give the morphism $B_t/SL_2
\longrightarrow Y_t/SL_2$ the structure of a vector bundle with
fiber $V^{\oplus k-3}$, i.e., of rank $2(k-3)$.  In particular $B_t/SL_2$
is $\mathbb{A}^1$-weakly equivalent to $Y_t/SL_2$ (e.g by Lemma
\ref{lem:contractibility}).

The condition that the hypersurface defined by $\varphi_t(\Delta_1, \Delta_2, \Delta_3) + 1 = 0$ be everywhere stable in $W$ is satisfied as long as the constant term is non-zero.  A generic hypersurface in $\mathbb{A}^3$ is smooth, and hence any $SL_2$-bundle $Y_t$ over it is smooth; it follows
that a generic choice of $\varphi_t$ determines a smooth everywhere
stable $B_t$ in $W$.  Furthermore, we can choose $\varphi_t$ so that the corresponding affine varieties occur in arbitrary dimensional families.  Indeed, we know that hypersurfaces in ${\mathbb P}^3$ of degree $d \geq 5$ are generically smooth and admit arbitrary dimensional moduli.  Fixing a hyperplane section $H \subset {\mathbb P}^3$, one can obtain hypersurfaces in ${\mathbb A}^3$ with the same property.

Let $Y_t/SL_2$ and $Y_s/SL_2$ be two such non-isomorphic
hyperbolic surfaces in $\mathbb{A}^3$.  If $B_t/SL_2$ is isomorphic
to $B_s/SL_2$ then $Y_t/SL_2$ and $Y_s/SL_2$ are $\aone$-weakly
equivalent; by assumption however, they are $\aone$-rigid and so would have to be isomorphic.  Thus $B_t/SL_2$ is not isomorphic to $B_s/SL_2$.
Consequently, any $p$-dimensional moduli of algebraically hyperbolic
surfaces in $\mathbb{A}^3$ induces $p$-dimensional moduli of
$B_t/SL_2$.
\end{proof}

\begin{lem}\label{lem:smoothfamilies}
There exist families as described in Theorem \ref{thm:moduli} such
that the quasi-affine quotients each admit a smooth affine closure
in $\Spec k[W]^{\ga}$ with smooth pure codimension $2$ boundary.
\end{lem}
\begin{proof}
In the proof of Lemma \ref{lem:families} the variety $B_t$ could be
chosen to be smooth, i.e., $F_{0,0} = 0$ is smooth. Now Lemma
\ref{lem:F00-smooth} implies $F = 0$ in $V \times W$ is smooth.  But
for $W = V^{\oplus k}$, there are no points with finite isotropy for
the $SL_2$ action on $V \times W$. Therefore $SL_2$ acts freely on
$F = 0$ in $V \times W$, so the quotient is smooth.
\end{proof}

\subsubsection*{Quotients of Dimension $\leq 2$}
\begin{claim}
\label{claim:dimension1}
There is a unique up to isomorphism $\aone$-contractible smooth scheme of dimension $1$; namely $\aone$.  Hence, any everywhere stable action of an $n$-dimensional unipotent group $U$ on ${\mathbb A}^{n+1}$ is isomorphic to $\aone$.
\end{claim}

\begin{proof}
Let $C$ be a smooth $\aone$-contractible curve.  Let $\bar{C}$ denote the projective completion of $C$ and let $g$ be the genus of this projective completion.  If the genus of $g$ is greater than $1$, then $C$ is $\aone$-rigid in the sense of \cite{MV} \S 2 Example 2.4, and hence not $\aone$-connected (thus, not $\aone$-contractible).  Therefore, $g$ must be $0$.

As $C$ is $\aone$-contractible, the set $C(k)$ is non-empty by \cite{MV} \S 3 Remark 2.5:  indeed, by \cite{MV} \S 2 Corollary 3.22, the canonical map $\bar{C}(L) \longrightarrow [\Spec L,\bar{C}]_{\aone}$ is surjective for any extension field $L/k$ (here $[\Spec L,X]_{\aone}$ denotes the set of $\aone$-homotopy classes of maps).  One knows that any smooth projective genus $0$ curve over a field $k$ possessing a $k$-rational point is isomorphic to $\pone$ over $k$.  Therefore, $C$ is a complement of a finite collection of $k$-rational points in $\pone$.  If $C$ is the complement of $\geq 2$ $k$-rational points in $\bar{C}$, then again $C$ is $\aone$-rigid.  Thus $C$ must be the complement in $\pone$ of a single $k$-rational point and is hence isomorphic to $\aone$.

Finally, by Corollary \ref{cor:main}, we know that the quotient ${\mathbb A}^{n+1}/U$ exists and is a smooth $\aone$-contractible scheme that is necessarily of dimension $1$.  Thus, by the previous paragraphs it must be isomorphic to $\aone$.
\end{proof}

\begin{claim}
\label{claim:dimension2}
Assume now that $k = \cplx$.  Suppose $U$ is an $n$-dimensional unipotent group.  The quotient of any everywhere stable action of $U$ on ${\mathbb A}^{n+2}$ is isomorphic to the affine plane.
\end{claim}

\begin{proof}
Suppose we have an everywhere stable action of $U$ on ${\mathbb A}^{n+2}$.  Let $X$ be the quotient ${\mathbb A}^{n+2}/U$.  Note that $X$ is a smooth quasi-affine $\aone$-contractible surface; in particular it is a smooth quasi-projective surface.  The topological space $X(\cplx)$ is therefore acyclic and hence by a result of Fujita (see \cite{Za} Lemma 2.1) $X$ is necessarily a smooth affine surface.  Furthermore, we know that the morphism ${\mathbb A}^{n+2} \longrightarrow X$ is in fact a trivial principal bundle isomorphic to $X_\cplx \times U$ by Corollary \ref{cor:main} (ii).  Now, since $X \times {\mathbb A}^{n} \cong {\mathbb A}^{n+2}$, it follows that $X$ is fact isomorphic to ${\mathbb A}^2$ by Fujita's proof of Zariski cancellation (see \cite{Fu}) in dimension $2$.
\end{proof}

\begin{rem}
If $X$ is a smooth algebraic variety, let $\bar{\kappa}(X)$ denote the logarithmic Kodaira dimension of $X$.  Suppose $f: Y \longrightarrow X$ is a $U$-torsor.  By the Iitaka Easy Addition theorem (see e.g. \cite{Za} Theorem 2.5 (c)), $\bar{\kappa}(Y) \leq \bar{\kappa}(U) + dim(X)$, and hence $\bar{\kappa}(Y)$ must be $-\infty$.  In other words, if a unipotent group acts everywhere stably on a variety $Y$, then $\bar{\kappa}(Y) = -\infty$.  However, the quotient of an everywhere stable action of a unipotent group on a smooth variety can have arbitrary logarithmic Kodaira dimension.  Indeed, let $X$ be a smooth affine variety with logarithmic Kodaira dimension $k$.  Then the usual translation action of $\ga$ on $\aone \times X$ (acting trivially on $X$) is everywhere stable with quotient isomorphic to $X$.
\end{rem}

\subsubsection*{Quotients of dimension $3$}
Again, suppose $k$ is the field $\cplx$.  In order to produce an example of a strictly quasi-affine quotient in dimension $3$ by our method, we would have to fix an embedding of ${\mathbb A}^4$ as a hypersurface in ${\mathbb A}^5$.  The restriction of the locally nilpotent derivation $D$ defining the $\ga$-action on ${\mathbb A}^5$ to ${\mathbb A}^4$ is triangular.  It follows from results of Deveney, Finston and van Rossum (see \cite{DFvR} Theorem 2.1) that {\em any} everywhere stable action of $\ga$ on ${\mathbb A}^4$ defined by a triangular locally nilpotent derivation has quotient isomorphic to ${\mathbb A}^3$.  This shows in particular, that the examples produced by the method of proof of Theorem \ref{thm:dimension4} are of minimal dimension.

\begin{question}
Does there exist a $3$-dimensional quasi-affine quotient of ${\mathbb A}^4$ by $\ga$?
\end{question}

\section{Consequences, Conjectures, and Comments}
\label{s:conjectures}
In this section, we emphasize some formal consequences of
$\aone$-contractibility and discuss some conjectures regarding the
structure of $\aone$-contractible smooth schemes.

\subsubsection*{Cohomology Computations}
The motivic homotopy category was constructed to study cohomology theories on the category of algebraic varieties.  For any $\aone$-contractible scheme $X$, and any space ${\mathcal Y}$, the sets of $\aone$-homotopy classes of maps $[X,{\mathcal Y}]_{\aone}$ and $[{\mathcal Y},X]_{\aone}$ are isomorphic to $[\Spec k,{\mathcal Y}]_{\aone}$ and $[{\mathcal Y},\Spec k]_{\aone}$ respectively.  Here are some trivial consequences of these facts.

 \begin{cor}
 Suppose $X$ is an $\aone$-contractible smooth scheme.  The structure map $X \longrightarrow \Spec k$ induces isomorphisms $H^{*,*}(\Spec k,\Z) \isomto H^{*,*}(X,\Z)$, $K^*(\Spec k) \isomto K^*(X)$, and $MGL^{*,*}(\Spec k) \isomto MGL^{*,*}(X)$.
 \end{cor}

\begin{proof}
All of these facts follow from the observation that the corresponding cohomology groups can be defined, unstably, as maps into an appropriate space: motivic Eilenberg-MacLane spaces $K(\Z(q),p)$ for motivic cohomology (see \cite{VICM}), $BGL_\infty$ for algebraic $K$-theory (see \cite{MV} \S4 Proposition 3.9), and $MGL_{\infty}$ for algebraic cobordism (see \cite{VICM}).
\end{proof}

 \begin{ex}
 Suppose $X$ is an $\aone$-contractible smooth scheme.  Then $Pic(X)$ is trivial, $\O^*(X)$ is isomorphic to $k^*$, and $K_2(X) \cong K_2^M(k)$ (where $K_i^{M}(k)$ is the $i$-th Milnor K-theory group of the field $k$).
 \end{ex}

\subsubsection*{Remarks about motives}
Suppose $X$ is any smooth algebraic variety.  Analogous to ordinary topology, one can study for $X$ motivic homology and motivic homology with compact supports (i.e. Borel-Moore motivic homology).  Voevodsky denotes the corresponding objects in the derived category of mixed motives by ${\sf M}(X)$ and ${\sf M}^c(X)$; these objects can be thought of as analogous to the usual singular chain complex and the singular chains with locally finite supports viewed as objects in the derived category of $\Z$-modules.  We refer the reader to \cite{VTriCat} for details about derived categories of motives.  We denote by $\dmeff{\Z}$ and ${\bf DM}_{gm}(\Spec k,\Z)$ Voevodsky's derived category of effective motivic complexes and derived category of geometric motives respectively (the Tate motive $\Z(1)$ is inverted in the latter category).

\subsubsection*{Consequences of a very general Hodge-type conjecture}
In this section we work with varieties over $\cplx$.  As we noted above, any smooth, $\aone$-contractible $\cplx$-algebraic variety has trivial motivic cohomology.  Via realization functors, this statement has consequences for the ordinary topology of such varieties.

The general form of the Hodge conjecture predicts that the embedding
of Grothendieck's category of homological motives (with
$\Q$-coefficients) into the category of $\Q$-Hodge structures is a
fully-faithful embedding. One can construct a Hodge realization
functor $R_{\mathcal H}$ on Voevodsky's derived category of motives
(see \cite{Hub} \S 3). Huber deduces the following result from
``standard conjectures" (see \cite{Hub} Proposition 3.4.1)

\begin{prop}
If $H^i(R_{\mathcal H}({\sf M}(X))) = 0$ for all $i$, then ${\sf M}(X)$ is equivalent to a point in the category ${\bf DM}_{gm}(\Spec k,\Q)$.
\end{prop}

In particular, this conjecture implies that if a smooth variety is
rationally acyclic, then it must also have rational motivic (co)homology
isomorphic to that of a point.  From this point of view, it is natural to study the motivic topology of varieties that are rationally acyclic.  In particular, the next question naturally presents itself.

\begin{question}
Do there exist topologically contractible smooth affine
$\cplx$-algebraic varieties that are not $\aone$-contractible?\footnote{In fact we can construct such
examples using affine modifications in the sense of Zariski; we
defer discussion to a future paper on affine contractible
varieties.}
\end{question}

\subsubsection*{Motivic topology at infinity: a dream}
As we noted in the introduction, open contractible $n$-manifolds $n \geq 4$ (PL or smooth) non-homeomorphic to $\real^n$ were necessarily non-simply connected at infinity.  However, rather than studying just the fundamental group at infinity, it is natural to study the whole ``homotopy type at infinity (see \cite{Ranicki} Chapter 9)."  A first step in this direction is to study the homology at infinity (see \cite{Ranicki} Chapter 3 for a definition of this notion).  It follows essentially from considerations involving Poincar\'e duality that the homology at infinity of a smooth open contractible $n$-manifold is that of an $n-1$-sphere.  The main goal of this section is to develop some notions of motivic topology at infinity.  The first thing to do is to study {\em motivic homology at $\infty$}; this notion has been introduced by Wildeshaus (see \cite{Wild}).

Notation as in \cite{VTriCat}, the motive ${\sf M}(X)$ and the compactly supported motive ${\sf M}^c(X)$ of a scheme $X$ are the objects in the derived category of motives $\dmeff{\Z}$ corresponding to the complexes $C_*\Z_{tr}(X)$ and $C_*\Z_{tr}^c(X)$.  There is a canonical morphism $\iota_X: \Z_{tr}(X) \longrightarrow \Z_{tr}^c(X)$.  The functor $C_*$ is exact and we can then define the {\em boundary motive} or the {\em motive at infinity}, denoted ${\sf M}^{\infty}(X)$ as the object in the derived category of motives corresponding to $C_*(Coker(\iota_X))[-1]$.

The motive at infinity is in some ways analogous to the {\em singular chain complex at infinity} (see \cite{Ranicki} Definition 3.8 (i)).  Similar to its topological analogue, the motive at infinity is one measure of the extent to which $X$ fails to be compact:  if $X$ is proper, the motive at infinity is trivial.  The following result shows that any $\aone$-contractible variety has motive at infinity that of a ``motivic sphere" of appropriate dimension.  Let ${\bf DM}(\Spec k,\Z)$ denote Voevodsky's derived category of mixed motives with the Tate motive inverted (see \cite{VTriCat} p. 192).

\begin{lem}
\label{lem:homatinfinity}
Suppose $X$ is an $m$-dimensional $\aone$-contractible finite type smooth scheme.  Then ${\sf M}^{\infty}(X) \cong \Z \oplus \Z(m)[2m-1]$ (non-canonically) as objects in ${\bf DM}(\Spec k,\Z)$.
\end{lem}

\begin{proof}
By definition of the motive at infinity, there is a distinguished triangle in ${\bf DM}(\Spec k,\Z)$ of the form
$$
{\sf M}^c(X)[-1] \longrightarrow {\sf M}^{\infty}(X) \longrightarrow {\sf M}(X) \longrightarrow {\sf M}^c(X).
$$
As $X$ is $\aone$-contractible, it follows that ${\sf M}(X) \cong \Z$.  Using motivic Poincar\'e duality (see \cite{VTriCat} Theorem 4.3.7), we see that
$${\sf M}^c(X) \cong {\sf M}(X)^{*} \tensor \Z(m)[2m] \cong \Z(m)[2m].$$
Finally, we know that $Hom_{{\bf DM}(\Spec k,\Z)}(\Z,\Z(m)[2m]) \cong H^{2m,m}(\Spec k,\Z)$ which vanishes for all $m > 0$.  Therefore, we can find a splitting of the distinguished triangle just mentioned and we obtain an isomorphism ${\sf M}^{\infty}(X) \cong \Z \oplus \Z(m)[2m-1]$.
\end{proof}

\begin{rem}
As the referee observed, this notion of motivic homology at infinity
is not sufficiently refined to distinguish, for example, phenomena
at infinity involving the real points.  Thus, the true
$\aone$-singular chain complex at infinity of a smooth variety
defined over $\real$ should at least see the topological singular
chain complex at infinity of the real points.  We refer the reader
to \cite{MorelFund} for a definition of the $\aone$-singular chain
complex of a variety $X$ and further development of the theory of
the $\aone$-fundamental group.  It would be interesting to formulate
an appropriate analogue of both of these objects ``at infinity."
\end{rem}

\section{Appendix}
\label{s:appendix}
In this appendix, we review some aspects of faithfully flat descent,
discuss Borel transfer, and prove Theorem \ref{thm:affinequotient}, which gives the cohomological criterion for quotients of affine varieties by unipotent group actions to be affine.  These results are used in the main body of the text but we have presented them here so as not to interrupt the narrative flow.

\subsubsection*{Faithfully flat descent}
Let us begin by recalling one form of faithfully flat descent; this material is well known, but we collect it here for the convenience of the reader.  A general reference for the material in this section is \cite{SGA1} Expose VIII.  Let $G$ be a group scheme.  For this section only, we will say $Y$ is a {\em left} (resp. right) $G$-scheme if $Y$ is a scheme equipped with a left (resp. right) $G$-action.  If we do not specify the chirality of an action in a theorem, we mean the result holds for both left and right actions.  We will use both left and right torsors.  Finally, if $Y$ is a $G$-scheme, we let ${\sf Qcoh}^G(Y)$ denote the category of $G$-equivariant sheaves on $Y$.

\begin{thm}
\label{thm:ffdescent}
Suppose $G$ is a linear algebraic group and that $f: {\mathscr P} \longrightarrow X$ is a $G$-torsor over a scheme $X$.  Then the functor
$$
f^*: {\sf Qcoh}(X) \longrightarrow {\sf Qcoh}^G({\mathscr P})
$$
is an equivalence of categories.  The functor $\F \mapsto (f_*(\F))^G$ defines an explicit quasi-inverse to $f^*$.
\end{thm}

\begin{proof}[Sketch of Proof]
Since $f: {\mathscr P} \longrightarrow X$ is a $G$-torsor, we have a canonical identification $G \times {\mathscr P} \isomto {\mathscr P} \times_X {\mathscr P}$.  Similarly, one obtains an isomorphism $G \times G \times {\mathscr P} \isomto {\mathscr P} \times_X {\mathscr P} \times_X {\mathscr P}$.  The morphism $f$ is faithfully flat, and one can check that specifying a descent datum on a quasi-coherent sheaf for $f$ is exactly the same as a specifying a $G$-equivariant structure on a quasi-coherent sheaf.  One then applies \cite{SGA1} Expose VIII Theorem 1.1.
\end{proof}

If $Y$ is a $G$-scheme, we let ${\mathscr A}^G(Y)$ denote the category of $G$-schemes affine over $Y$ such that the structure morphism is $G$-equivariant.  Suppose now that $G$ is a linear algebraic group and $H$ is a closed algebraic subgroup.  Then we know the homogeneous space quotient $G/H$ exists and that the morphism $\pi: G \longrightarrow G/H$ is a right $H$-torsor which is $G$-equivariant for the natural left $G$-actions on $G$ and $H$.  Furthermore, the structure morphism $s:G \longrightarrow \Spec k$ is a left $G$-torsor.

\begin{cor}
\label{cor:sheafproperty}
The functor $\F \mapsto f_*(s^*\F)^H$ defines an equivalence of categories $${\sf Qcoh}^H(\Spec k) \isomto {\sf Qcoh}^G(G/H).$$  This functor also determines an equivalence of categories $${\mathscr A}^H(\Spec k) \isomto {\mathscr A}^G(G/H).$$An inverse to the last functor can be obtained by taking the scheme-theoretic fiber of a morphism $f: X \longrightarrow G/H$ over the identity coset of $G/H$.
\end{cor}

\begin{proof}
The first statement follows immediately from Theorem \ref{thm:ffdescent}.
The second statement follows because there is an equivalence of categories between schemes affine over $Y$ and quasi-coherent sheaves of $\O_Y$-algebras.  For more details see \cite{SGA1} Expose VIII Thm 2.1.
\end{proof}

For an affine $H$-scheme $X$, we write $G *_H X$ for the image of $X$ under the functor of Corollary \ref{cor:sheafproperty}; the scheme $G *_H X$ is referred to as a contracted product scheme.  By construction it can be identified with the quotient of the product $G \times X$ by the action of $H$ defined by $h \cdot (g,x) = (gh,h \cdot x)$.  Given a morphism $f: X \longrightarrow X'$ of affine $H$-schemes, we denote by $G *_H f$ the induced morphism $G *_H X \longrightarrow G *_H X'$.

\begin{rem}
\label{rem:quasi-affine}
In fact by Corollaire \cite{SGA1} Expose VIII Cor. 7.9, if $X$ is a {\em quasi-affine} $H$-scheme, then $G *_H X$ exists as a scheme.
\end{rem}

\begin{rem}
Through the paper we have assumed that the term ``scheme" is synonymous with separated scheme, locally of finite type; for this remark we lift this convention.  Theorem \ref{thm:ffdescent} and Corollary \ref{cor:sheafproperty} above make no reference to any finiteness assumptions and hold as well for affine schemes that are not locally of finite type.
\end{rem}

Let $P(f)$ be a property of a morphism $f: X \longrightarrow Y$ of schemes that is stable by base change and local for the \'etale topology.  The following properties of morphisms of schemes are of this form:
\begin{itemize}
\item surjective, radiciel, universally bijective, universally open, universally submersive, separated, quasi-compact, locally of finite type, finite type, open immersion, closed immersion, affine, quasi-affine, integral, geometrically irreducible fibers, geometrically reduced, geometrically connected fibers, flat, or smooth.
\end{itemize}
For a more complete list of such properties, we refer the reader to \cite{SGA1} Expose VIII \S 3-8.

\begin{cor}
\label{cor:morphismproperty}
A morphism of quasi-affine $H$-schemes $f: X \longrightarrow X'$ has a property $P(f)$ as in the previous paragraph if and only if the induced morphism $G *_H f$ has the same property.
\end{cor}

\begin{proof}
This follows from Corollary \ref{cor:sheafproperty} and Remark \ref{rem:quasi-affine}.  We refer the reader to \cite{SGA1} Exposes VIII \S 3, \S 4.
\end{proof}

\begin{rem}
\label{rem:qaffmorphisms}
The construction of Remark \ref{rem:quasi-affine} and Corollary \ref{cor:morphismproperty} then gives an equivalence of categories between the category of $G$-schemes quasi-affine over $G/H$ (with quasi-affine morphisms) and the category of quasi-affine $H$-schemes (with quasi-affine morphisms).
\end{rem}

\subsubsection*{Zariski local triviality of $U$-torsors, proof of Corollary \ref{cor:unipgrptriv}}
\label{ss:groththeorem}
By our definitions, $U$-torsors on a scheme $X$ are classified by the flat cohomology group $H^1_{fppf}(X,U)$.  Observe that by definition $H^1_{fppf}(X,\ga) \cong H^1_{fppf}(X,\O_X)$.  By Theorem \ref{thm:flag}, $U$ admits an increasing filtration by algebraic subgroups with subquotients isomorphic to $\ga$.  Thus the sheaf defined by $U$ on $X$ is a coherent sheaf (being a successive extension of sheaves isomorphic to $\O_X$).  The result then follows from the fact that for any coherent sheaf $\F$, the canonical map $H^1_{Zar}(X,\F) \longrightarrow H^1_{fppf}(X,\F)$ is an isomorphism (see e.g. \cite{Milne} Chapter III Proposition 3.7).

\subsubsection*{Naturality of contracted products}
We now explore the naturality properties of the contracted product construction of the previous section.  In the main body of the text, we will only consider the situation where $H$ is a connected unipotent group or $H$ is a reductive group.

\begin{thm}
\label{thm:matsushima}
Let $G$ be a reductive linear algebraic group.  Then the quotient $G/H$ is affine if and only if $H$ is a reductive subgroup of $G$.  If $U$ is a unipotent subgroup of $G$, then $G/U$ is strictly quasi-affine.
\end{thm}

\begin{proof}[Sketch of Proof]
The first part of this result is the statement of Matsushima's theorem (see \cite{Hab} Theorem 3.3).  Matsushima's theorem also guarantees that if $U$ is unipotent, then $G/U$ is not affine.  One can check that for any $k$-defined parabolic $P$, with unipotent radical $R_u(P)$, the quotient $G/R_u(P)$ is strictly quasi-affine.  Since $U$ may be embedded in $R_u(P)$ for some $P$, and the quotient $R_u(P)$ is isomorphic to affine space, we get an \'etale locally trivial fiber bundle $G/U \longrightarrow G/R_u(P)$ with fibers isomorphic to affine space.  Since the base of this fibration is
quasi-affine, it follows that $G/U$ is necessarily quasi-affine by \ref{cor:morphismproperty}.  It follows that $G/U$ is strictly quasi-affine.
\end{proof}

Suppose we have a sequence of inclusions $U \hookrightarrow G \hookrightarrow G'$ where $G'$ and $G$ are reductive and $U$ is quasi-affine.  If $X$ is a $U$-scheme, then it follows immediately from Corollary \ref{cor:morphismproperty} and \ref{thm:matsushima} that both $G *_U X$  and $G' *_U X$ are necessarily quasi-affine schemes.  Indeed, the morphism $X \longrightarrow \Spec k$ is affine, and the composite of an affine morphism and quasi-affine morphism is quasi-affine.  Furthermore we have a canonical identification $G' *_G (G *_U X) \isomto G' *_U X$.   We now study global sections of the structure sheaf on a contracted product scheme.

\subsubsection*{Borel Transfer}
\label{ss:transfer}
Henceforth let us write $k[Y]$ for the $k$-algebra of global sections $\Gamma(Y,\O_Y)$ for any quasi-affine scheme $Y$.  Suppose $X$ is a quasi-affine $U$-scheme.  The closed immersion $X \hookrightarrow G \times X$ defined by $x
\mapsto (e,x)$ descends to a closed immersion $\iota: X
\hookrightarrow G *_U X$.  Pull-back by $\iota$ defines a morphism $k[G *_U X] \longrightarrow
k[X]$ sending $G$-invariant functions to $U$-invariant functions.
The Borel transfer principle says that induced map $$\iota^*: k[G
*_U X]^G \isomto k[X]^U$$ is an isomorphism (this follows from e.g. \cite{SGA1} Expose VIII Cor. 1.7).  This isomorphism is
functorial in the following sense.  Suppose $G'$ is another
reductive group such that $G \subset G'$.  We then have induced
morphisms $\iota': X \hookrightarrow G'*_U X$ and $\iota_1: G*_U X
\hookrightarrow G' *_G (G *_U X)$.  Finally, we have
a canonical identification $\iota_1 \circ \iota = \iota'$. Through
the rest of this document, we will use the induced isomorphisms
without comment.

\subsubsection*{Proof of Theorem \ref{thm:affinequotient}}
\label{ss:affinequotient}
Let us now prove Theorem \ref{thm:affinequotient}.  Essentially, this theorem is a long exercise in faithfully flat descent and we felt this justified its deferral to this appendix.
Suppose $\pi: X \longrightarrow X/U$ is a principal bundle (in particular it is faithfully flat and affine).  We can consider the \u Cech simplicial scheme $\breve{C}(\pi)$ attached to this morphism:
\begin{equation}
\cdots X \times_{X/U} X \rightrightarrows X \longrightarrow X/U
\end{equation}
In other words, the $n^{th}$-term of $\breve{C}(\pi)$ is the $n+1$-fold fiber product of $X$ over $X/U$; the partial projections and the relative diagonals give the face and degeneracy maps.  This is a simplicial scheme augmented toward $X/U$.  Since $X \longrightarrow X/U$ is a $U$-torsor, the \u Cech simplicial scheme is isomorphic as a simplicial scheme to the bar simplicial scheme whose $n^{th}$ term is $X \times U^{\times n}$.  We denote this scheme by $X \times U^{\bullet}$.  As $X$ and $\pi$ are affine, it follows that $\breve{C}(\pi)$ is a simplicial affine scheme.

Given a quasi-coherent sheaf $\F$ on $X/U$, we can consider the simplicial sheaf on $\F$ on $\breve{C}(\pi)$ induced by pull-back.  As $\breve{C}(\pi)$ is an affine simplicial scheme, the spectral sequence of cohomological descent attached to the morphism $\pi$ (see \cite{SGA42} V.bis 2.5) degenerates and gives us a complex
\begin{equation}
\breve{C}(\pi,\F) = \Gamma(X/U,\F) \longrightarrow \Gamma(X,\pi^*\F) \longrightarrow \Gamma(X \times U, p_1^*\F) \longrightarrow \cdots
\end{equation}
whose $i$-th cohomology computes the group $H^i(X,\F)$.  Furthermore, the isomorphism of simplicial schemes described in the previous paragraph identifies the groups $H^i(X/U,\F)$ with the group cohomology $H^i(U,\Gamma(X,\pi^*\F))$.  With these notations, our proof is fairly streamlined.

\begin{proof}[Proof of Theorem \ref{thm:affinequotient}]
(i $\Longrightarrow$ ii).  By the identifications of the previous paragraph, we see that $0 = H^1(X/U,\O_{X/U}) = H^1(U,\Gamma(X,\O_X))$.  Now we use the van Est spectral sequence; let us recall the setup (see \cite{Teleman} \S 6.1).  Let ${\sf Rep}(U)$ denote the category of locally finite algebraic $U$-modules, and similarly for ${\sf Rep}({\mathfrak u})$.  There is a restriction functor $Res^U_{{\mathfrak u}}: {\sf Rep}(U) \longrightarrow {\sf Rep}({\mathfrak u})$, and this functor possesses a right adjoint functor $Ind_{\mathfrak u}^U$ of induction.  Let $(\cdot)^U$ denote the functor of $U$-invariants and let $(\cdot)^{\mathfrak u}$ denote the functor of ${\mathfrak u}$-invariants.  The composite functor $(Ind_{\mathfrak u}^U(\cdot))^U$ is isomorphic to the functor $(\cdot)^{\mathfrak u}$.  If $V$ is a ${\mathfrak u}$-module, we obtain a Grothendieck spectral sequence:
\begin{equation} E^{p,q}_2 = H^p(U,{\bf R}^q Ind(V))
\Longrightarrow H^*({\mathfrak u},V).
\end{equation}
The $U$-action can be untwisted, and we obtain an isomorphism
\begin{equation}
E^{p,q}_2 \cong H^p(U,H^q_{dR}(U) \tensor V) \Longrightarrow H^*({\mathfrak u},V).
\end{equation}
where $H^q_{dR}(U)$ is the algebraic de Rham cohomology of $U$.  Since $U$ is connected, it acts trivially on $H_{dR}^*(U)$ and each $E_2^{p,q}$ factors as a tensor product: $H^p(U, H^q(U) \tensor V) \cong H^p_{dR}(U) \tensor H^p(U,V)$.  The algebraic de Rham cohomology of $U$ is trivial by homotopy invariance: $U$ is isomorphic to affine space and hence the spectral sequence degenerates and defines an isomorphism of ${\mathfrak u}$-cohomology and $U$-cohomology.  In particular, we see that $H^1(U,\Gamma(X,\O_X))$ vanishes if and only if $H^1({\mathfrak u},\Gamma(X,\O_X))$-vanishes.

(ii) $\Longrightarrow$ (iii).  We can't assume that the quotient $X/U$ exists as a scheme, but as $U$ is a smooth $k$-group-scheme and $X$ is a $k$-variety we can consider the Artin quotient stack $[X/U]$ (see \cite{LMB} Chapter 3).  The stack $[X/U]$ comes equipped with the universal atlas $X \longrightarrow [X/U]$.  Furthermore, the coherent cohomology of any sheaf on $[X/U]$ can be computed as the cohomology of the associated \u Cech complex as described in the paragraph just preceding this proof (See e.g. \cite{LMB} \S 13.5).  Using the argument of the previous paragraph identifying the ${\mathfrak u}$-cohomology of $\Gamma(X,\O_X)$ and the $U$-cohomology of $\Gamma(X,\O_X)$, we see that $H^1([X/U],\O_{[X/U]})$ vanishes.

Since $H^1([X/U],\O_{[X/U]})$ vanishes, all $U$-torsors on $[X/U]$ must be trivial, and hence the morphism $\pi: X \longrightarrow [X/U]$ admits a section, i.e. a morphism $s: [X/U] \longrightarrow X$.  We claim that such a morphism is necessarily representable.  Indeed, as $s$ is a section, $s$ has the property that the two composite morphisms are representable morphisms: $s\pi = Id_{[X/U]}$ and $\pi s$ is a morphism of schemes.  The morphism $X \times_{[X/U]} X \longrightarrow X \times X$ is representable as well:  it is canonically isomorphic to the action map $X \times U \longrightarrow X \times X$.  It follows by \cite{LMB} Lemme 3.12 c(i), the section $s$ is a representable morphism.  Now, any stack admitting a representable morphism to a scheme must actually be an algebraic space by \cite{LMB} Lemme 3.12(a).  It suffices to check that the algebraic space $[X/U]$ is in fact an affine scheme.  In this case, observe that the morphism $X \longrightarrow X \times U$ is an affine morphism and thus by descent, the morphism $[X/U] \longrightarrow X$ is an affine morphism.  Since $X$ is an affine scheme it follows by \cite{Knutson} Chapter 2 Extension 3.8 that $[X/U]$ is necessarily an affine scheme.

(iii) $\Longrightarrow$ (i).  Suppose $\pi: X \longrightarrow X/U$ is a trivial $U$-torsor.  In this case, we have an isomorphism $X \isomto X/U \times U$.  The unit morphism $\Spec k \longrightarrow U$ gives rise to a map $X/U \longrightarrow X$ that is a closed immersion and thus $X/U$ is an affine scheme.
\end{proof}

\begin{rem}
As a corollary of our proof, we see that the group $H^1({\mathfrak u},k[X])$ is an {\em invariant} of the action because it is identified with the group $H^1([X/U],U)$.
\end{rem}

\begin{footnotesize}
\bibliographystyle{alpha}
\bibliography{exoticaffines}
\end{footnotesize}

\end{document}